\documentclass[a4paper,12pt]{article}
\usepackage[cp1251]{inputenc} 
\usepackage[english]{babel}
\usepackage{amsfonts,amssymb}
\title{Differential modules with
$\infty$-simplicial faces and $A_\infty$-algebras}
\author{S.V. Lapin}
\date{}
\newcommand{\F}{F_{\infty}}
\newcommand{\D}{D_{\infty}}
\newcommand{\bu}{\bullet}
\newcommand{\E}{E_{\infty}}
\newcommand{\A}{A_{\infty}}

\newcommand{\Hom}{{\rm Hom}}
\newcommand{\p}{\partial}

\newcommand{\sign}{\rm sign}
\begin{document}
\maketitle

\begin{abstract} In the present paper, by using the colored version of the Koszul duality, the concept
of a differential module with $\infty$-simplicial faces is
introduced. The homotopy invariance of the structure of a
differential module with $\infty$-simplicial faces is proved. The
relationships between differential modules with
$\infty$-simplicial faces and $A_\infty$-algebras are established.
The notion of a chain realization of a differential module with
$\infty$-simplicial faces and the concept of a tensor product of
differential modules with $\infty$-simplicial faces are
introduced. It is proved that for an arbitrary $A_\infty$-algebra
the chain realization of the tensor differential module with
$\infty$-simplicial faces, which corresponds to this
$A_\infty$-algebra, and the $B$-construction of this
$A_\infty$-algebra are isomorphic differential coalgebras.

\end{abstract}

The concept of a chain realization of simplicial differential
modules and in particular the concept of a chain realization of
differential modules with simplicial faces was introduced in
\cite{Smir0}, when studying homotopy properties of
$\E$-coalgebras. As was showed in [1], a chain realization of the
tensor differential module with simplicial faces, which is defined
by an arbitrary differential algebra, is a $B$-construction of
this differential algebra. By using the representation of a
$B$-construction as a chain realization, in \cite{Smirnov} was
obtained that the $B$-construction of an arbitrary cocommutative
Hopf algebra is an $\E$-algebra, which by a suitable "gluing" of
$\E$-coalgebras of chain complexes of standard simplicial
simplexes is obtained. The presence of the structure
$\E$-coalgebra on the above considered $B$-construction implies
that on cohomology this $B$-construction exist the natural action
of Steenrod operations, which satisfy the Adem relations
\cite{Smir0} and Cartan formula \cite{Lap0}. In particular, if we
consider the Steenrod algebra as a cocommutative Hopf algebra,
then we obtain the well known action of Steenrod operations on
cohomology of the Steenrod algebra (see for example \cite{May}),
i.e., on the second term of the Adams spectral sequence for the
stable homotopy groups of spheres.

The present paper consists of three paragraphs. In the first
paragraph, the necessary notions, constructions and assertions of
the Koszul duality theory for quadratic algebras \cite{LV} and the
differential Lie modules theory \cite{Smir1} respectively are
carried over to the cases of quadratic colored algebras and
colored graded coalgebras. In the second paragraph, a description
of the colored coalgebra $(F^!,\nabla)$ Koszul dual to the
quadratic colored algebra of simplicial faces $(F,\pi)$ is given.
Further, the notion of a differential module with
$\infty$-simplicial faces is introduced as a differential Lie
module over a colored coalgebra $(F^!,\nabla)$ or, equivalently,
as a differential module over a co-$B$-construction of this
colored coalgebra. The homotopy invariance of the structure of a
differential module with $\infty$-simplicial faces is established.
The connection between the concepts of a differential module with
$\infty$-simplicial faces and differential module with homotopy
simplicial faces \cite{Lap3} is considered. In the third
paragraph, the concepts of a chain realization of a differential
module with $\infty$-simplicial faces and tensor product of
differential modules with $\infty$-simplicial faces are
introduced. It is shown that every $\A$-algebra defines the tensor
differential module with $\infty$-simplicial faces. Further, for
an arbitrary $\A$-algebra, it is proved that the chain realization
of its tensor differential module with $\infty$-simplicial faces
and the $B$-construction of this $\A$-algebra are isomorphic as
differential coalgebras.

We proceed to precise definitions and statements. All modules and
maps of modules considered in this paper are assumed to be
$K$-modules and $K$-linear maps of modules, respectively, where
$K$ is an arbitrary commutative ring with unity. \vspace{0.5cm}

\centerline{\bf \S\,1. Colored coalgebras Koszul dual to quadratic
colored algebras} \vspace{0.5cm}

Let $I$ be a set of nonnegative integers. Each element of the set
$I$ we will called a color. A colored differential module with
colors from the set $I$ or, more briefly, simply a colored
differential module $(X,d)$ is by definition an arbitrary family
of differential graded modules $X=\{X(s,t)_m\}$, $m\in\mathbb{Z}$,
$d:X(s,t)_\bullet\to X(s,t)_{\bullet-1}$ that are indexed by all
pairs of elements $(s,t)\in I\times I$. In the case when colored
differential module $(X,d)$ satisfies the condition $d=0$ for all
$s,t\in I$, we will say that $X$ is a colored graded module.

A map of colored differential modules $f:(X,d)\to (Y,d)$ is any
family of maps of differential graded modules
$f=\{f(s,t):(X(s,t),d)\to (Y(s,t),d)\}_{s,t\in I}$. Similarly we
introduce the notion of a homotopy $h:X_\bu\to Y_{\bu+1}$ between
maps colored differential modules $f,g:(X,d)\to (Y,d)$. It is
clearly that colored differential modules and their maps form a
category in which obvious way direct sums and direct products of
objects, the kernels and cokernels of morphisms, and also homology
of objects are defined.

In further a differential bigraded module we will call an
arbitrary differential bigraded module $(X,d)$, $X=\{X_{n,m}\}$,
$n,m\in\mathbb{Z}$, $n\geqslant 0$, $d:X_{*,\bu}\to X_{*,\bu-1}$.
For an arbitrary differential bigraded modules $(X,d)$ and $(Y,d)$
there is the colored differential module $(\hom(X;Y),d)$. The
elements of the module $\hom(X;Y)(s,t)_m$ are an arbitrary maps of
graded modules $f:X_{t,\bullet}\to Y_{s,\bullet+m}$ of degree
$m\in\mathbb{Z}$, and the differential
$d:\hom(X;Y)(s,t)_\bullet\to\hom(X;Y)(s,t)_{\bullet-1}$ for any
fixed $s,t\in I$ is given on elements $f\in \hom(X;Y)(s,t)_m$ by
the following formula: $$d(f)=df+(-1)^{s-t+m+1}fd:X_{t,\bullet}\to
Y_{s,\bullet+m-1}.$$

For an arbitrary colored differential modules $(X,d)$ and $(Y,d)$,
we set by definition $X(s,k)_p\otimes Y(l,t)_q=0$, $k\not=l$, and
define $X(s,k)_p\otimes Y(l,t)_q$, $k=l$,  as the usual tensor
product of modules.

A tensor product of the colored differential modules $(X,d)$ and
$(Y,d)$ is the colored differential module $(X\otimes Y,d)$, where
$$(X\otimes Y)(s,t)_m=\bigoplus\limits_{k,l\in
I}\bigoplus\limits_{p+q=m}X(s,k)_p\otimes
Y(l,t)_q=\bigoplus\limits_{k\in
I}\bigoplus\limits_{p+q=m}X(s,k)_p\otimes Y(k,t)_q,$$ and its
differential $d:(X\otimes Y)(s,t)_\bu\to (X\otimes
Y)(s,t)_{\bu-1}$, for an arbitrary element $x\otimes y\in
X(s,l)_p\otimes Y(l,t)_q$, is defined by the following formula:
$$d(x\otimes y)=d(x)\otimes y+(-1)^{s-l+p}x\otimes d(y).$$

A colored differential algebra $(A,d,\pi)$ is an arbitrary colored
differential module $(A,d)$ equipped with a multiplication
$\pi:A\otimes A\to A$, which is a map of colored differential
modules and satisfies the associativity condition $\pi(\pi\otimes
1)=\pi(1\otimes\pi)$. If a colored differential algebra
$(A,d,\pi)$ satisfies the condition $d=0$ then we will say that
there is given a colored graded algebra $(A,\pi)$.

A map of colored differential algebras $f:(A',d,\pi)\to
(A'',d,\pi)$ is an arbitrary map of colored differential modules
$f:(A',d)\to (A'',d)$ that satisfies the condition $\pi(f\otimes
f)=f\pi$.

A unity of a colored differential algebra $(A,d,\pi)$ is a family
of elements $1_*=\{1_k\}$, $k\in I$, $1_k\in A(k,k)_0$, that
satisfy the condition $\pi(1_s\otimes a)=a=\pi(a\otimes 1_t)$ for
each element $a\in A(s,t)_m$, $s,t\in I$, $m\in\mathbb{Z}$.

The canonical example of a colored differential algebra is the
colored differential algebra $(\hom(X;X),d,\pi)$, where
$(\hom(X;X),d)$ is considered above colored differential module,
which for any differential bigraded module $(X,d)$ is defined. A
multiplication $\pi:\hom(X;X)\otimes\hom(X;X)\to\hom(X;X)$ of this
colored differential algebra is given for an arbitrary elements
$g\in\hom(X,X)(s,k)_n$ and $f\in\hom(X,X)(l,t)_m$ by the following
rule: $$\pi(g\otimes f)=\left\{\begin{array}{ll}
gf\in\hom(X;X)(s,t)_{n+m},&k=l,\\
0\in\hom(X;X)(s,t)_{n+m},&k\not=l,\\
\end{array}\right.$$
where $gf:X_{t,\bullet}\to X_{s,\bullet+n+m}$ is a composition of
the maps $f:X_{t,\bullet}\to X_{(l=k),\bullet+m}$ and
$g:X_{(k=l),\bullet}\to X_{s,\bullet+n}$. It is clear that a unity
of the colored differential algebra $(\hom(X,X),d,\pi)$ is a
family of elements $1_*=\{1_s\}$, $s\in I$, where for any $s\in I$
the element $1_s\in\hom(X;X)(s,s)_0$ is the identity map of the
graded module $X_{s,\bullet}$.

Further we will denote by $K_I$ the colored graded module that is
defined by equalities $K_I(s,s)_m=K$, $m=0$, $s\in I$,
$K_I(s,s)_m=0$, $m\not=0$, $s\in I$, and $K_I(s,t)_m=0$,
$s\not=t$, $m\in\mathbb{Z}$. It is clearly that the multiplication
in the ring $K$ allows regarded $K_I$  as colored graded algebra
$(K_I,\pi)$.

By using considered above the operation of a tensor product of
colored graded modules, we define the usual way the notions of a
tensor algebra of a graded module, a two-sided ideal of a colored
graded algebra generated by a submodule of this colored graded
algebra and a quotient algebra of a colored graded algebra by a
two-sided ideal.

{\bf Definition 1.1}. A colored graded algebra $(A,\pi)$ with a
unity is called a quadratic colored graded algebra or, more
briefly, a quadratic colored algebra if this colored graded
algebra $(A,\pi)$ is isomorphic to the quotient algebra $R =
T(M)/(Q)$, where $T(M)$ is a tensor algebra of some colored graded
module $M=\{M{(s,t)_m}\}_{s,t\in I}$, $m\in\mathbb{Z}$,
$m\geqslant 0$, and $(Q)$ is a two-sided ideal of the colored
graded algebra $T(M)$ generated by some submodule $Q$ of the
colored graded module $M\otimes M\subset T(M)$.

We note that for any quadratic colored algebra $R=T(M)/(Q)$ the
composition of the obvious embedding $M\to T(M)$ and the
projection $T(M)\to T(M)/(Q)$ is an embedding $M\to R=T(M)/(Q)$ of
colored graded modules. In what follows, we identify the image of
this embedding with the colored graded module $M$ and always
assume that $M$ is a submodule of the colored graded module $R$.

A colored graded coalgebra $(C,\nabla)$ is an arbitrary colored
graded module $C=\{C(s,t)_m\}_{s,t\in I}$, $m\in\mathbb{Z}$,
$m\geqslant 0$, equipped with a comultiplication $\nabla:C\to
C\otimes C$ that is a map of colored graded modules satisfies the
condition $(\nabla\otimes 1)\nabla=(1\otimes\nabla)\nabla$. A
colored graded coalgebra $(C,\nabla)$ is said to be connected if
$C_0=(K_I)_0$. For a colored graded coalgebra $(C,\nabla)$ the
usual way we define notions of counit $\varepsilon:C\to K_I$ and
coaugmentation $\nu:K_I\to C$.

It is easy to see that a comultiplication in the colored graded
coalgebra $(C,\nabla)$ with counit $\varepsilon:C\to K_I$ and
coaugmentation determines a comultiplication
$\overline{\nabla}:\overline{C}\to\overline{C}\otimes\overline{C}$
in the colored graded submodule $\overline{C}=\ker(\varepsilon)$,
which for an arbitrary element $c\in\overline{C}(s,t)$, $s,t\in
I$, is defined by the formula
$\overline{\nabla}(c)=\nabla(c)-(1_s\otimes c+c\otimes 1_t)$.
Thus, for any colored graded coalgebra $(C,\nabla)$ the colored
graded coalgebra $(\overline{C},\overline{\nabla})$ without counit
is always defined.

The suspension of a colored graded module $M$ is a colored graded
module $SM$ with grading defined for $s,t\in I$ by the equality
$(SM)(s,t)_{m+1}=M(s,t)_{m}$. Let us denote the elements of $SM$
by $[x]$, where $x\in M$. For any submodule $Q$ of the colored
graded module $M\otimes M$ we denote by $S^{\otimes 2}Q$ the
corresponding submodule of the colored graded module $SM\otimes
SM$.

{\bf Definition 1.2}. The colored coalgebra Koszul dual to a
quadratic colored algebra $R=T(M)/(Q)$ is the colored graded
coalgebra $(R^!,\nabla)$ with counit $\varepsilon:R^!\to K_I$ and
coaugmentation $\nu:K_I\to R^!$ defined by the following formulae:
$$R^!=\bigoplus_{k\geqslant
0}(R^!)^{(k)},~~(R^!)^{(0)}=K_I,~~(R^!)^{(1)}=SM,$$
$$(R^!)^{(k)}=\bigcap_{i+2+j=k}(SM)^{\otimes i}\otimes S^{\otimes
2}Q\otimes (SM)^{\otimes j},~~k\geqslant 2,~i\geqslant
0,~j\geqslant 0,$$ $$\nabla(1_s)=1_s\otimes 1_s,~~~1_s\in
(R^!)^{(0)}(s,s)_0=K_I(s,s)_0,~s\in I,$$
$$\nabla([x_1\otimes\dots\otimes x_k])=1_s\otimes
[x_1,\dots,x_k]+[x_1,\dots,x_k]\otimes 1_t+$$
$$+\sum_{i=1}^{k-1}\,[x_1,\dots,x_i]\otimes
[x_{i+1},\dots,x_k],~~~[x_1,\dots,x_k]\in (R^!)^{(k)}(s,t),~s,t\in
I,$$
$$\varepsilon(1_s)=1_s,~\varepsilon([x_1,\dots,x_k])=0,~k\geqslant
1,~~\nu(1_s)=1_s,$$ where
$[x_1,\dots,x_k]=[x_1]\otimes\dots\otimes [x_k]$, $x_j\in
F^!(s_j,t_j)$, $s_j,t_j\in I$, $1\leqslant j\leqslant k$, $s_1=s$,
$t_k=t$, $s_p=t_{p+1}$, $1\leqslant p\leqslant k-1$.

Let us now consider the notion of a co-$B$-construction of a
connected colored graded coalgebra with counit and coaugmentation.

A desuspension over a colored graded module $M$ is the colored
graded module $S^{-1}M$, which for any $s,t\in I$ is defined by
the formula $(S^{-1}X)(s,t)_{m-1}=X(s,t)_m$. The elements of
$S^{-1}X$ are traditionally denoted by $[x]$, where $x\in M$.

Let $(\overline{C},\overline{\nabla})$ be the colored graded
coalgebra without counit determined by the colored graded
coalgebra $(C,\nabla)$ with counit and coaugmentation. A
co-$B$-construction of a colored graded coalgebra $(C,\nabla)$ is
the colored differential algebra $(\Omega(C),d,\pi)$ defined as
follows: $$\Omega(C)=\bigoplus_{k\geqslant
0}(S^{-1}\overline{C})^{\otimes k},~~(S^{-1}\overline{C})^{\otimes
0}=K_I,$$
$$d([c_1,\dots,c_k])=\sum_{i=1}^{k}(-1)^{i+\mu_i}[c_1,\dots,c_{i-1},
\overline{\nabla}(c_i),c_{i+1},\dots,c_k],$$
$$\overline{\nabla}(c_i)=\sum c_i'\otimes c_i'',\quad
\pi([c_1,\dots,c_q]\otimes
[c_{q+1},\dots,c_k])=[c_1,\dots,c_q,c_{q+1},\dots,c_k],$$ where
$[c_1,\dots,c_k]=[c_1]\otimes\dots\otimes[c_k]$,
$c_j\in\overline{C}(s_j,t_j)_{q_j}$, $s_j=t_{j+1}$, $1\leqslant
j\leqslant k$, $c_i'\in\overline{C}(s_i,t_i')_{q_i'}$,
$\mu_i=s_1-t_i'+q_1+\dots+q_{i-1}+q_i'$, $1\leqslant i\leqslant
k$.

A twisting cochain from a colored graded coalgebra $(C,\nabla)$ to
a colored differential algebra $(A,d,\pi)$ is defined as a map
$\varphi:C_\bu\to A_{\bu-1}$ of colored graded modules which has
degree $-1$ and satisfies the cochain twisting condition
$d\varphi+\varphi\cup\varphi=0$, where the map of colored graded
modules $\varphi\cup\varphi:C_\bu\to A_{\bu-2}$ of degree $-2$ is
defined by the formula
$\varphi\cup\varphi=\pi(\varphi\otimes\varphi)\nabla$.

The simplest example of a twisting cochain from any colored graded
coalgebra $(C,\nabla)$ to the colored differential algebra
$(\Omega(C),d,\pi)$ is the map $\varphi^{\Omega}:C_\bu\to
\Omega(C)_{\bu-1}$, $\varphi^{\Omega}(c)=[c]$, $c\in C$.

For any twisting cochain $\varphi:C_\bu\to A_{\bu-1}$ from a
colored graded coalgebra $(C,\nabla)$ to a colored graded algebra
$(A,\pi)$ is defined the map of colored differential algebras
$\Omega(\varphi):(\Omega(C),d,\pi)\to (A,d=0,\pi)$ that on the
generators $[c_1,\dots,c_k]$ of the colored graded module
$\Omega(C)$ is given by the following formula:
$$\Omega(\varphi)([c_1,\dots,c_k])=\pi^{(k)}(\varphi(c_1)\otimes\dots\otimes\varphi(c_k)),$$
where $\pi^{(1)}$ is the identity map of the colored algebra $A$,
and for $k\geqslant 2$ the map
$\pi^{(k)}=\pi(1\otimes\pi^{(k-1)})$ is the iterated
multiplication in the colored algebra $(A,\pi)$.

Let us now consider the notion of a differential module over a
colored differential algebra.

A tensor product of a colored differential module $(X,d)$ and
differential bigraded module $(Y,d)$ is the differential bigraded
module $(X\otimes Y,d)$ for which $$(X\otimes
Y)_{n,m}=\bigoplus\limits_{s\in I\atop{p+q=m}}X(n,s)_p\otimes
Y_{s,q},~~n\in I,~m\in\mathbb{Z},$$ and the differential
$d:(X\otimes Y)_{n,\bullet}\to (X\otimes Y)_{n,\bullet-1}$ is
defined on an arbitrary element $x\otimes y\in X(n,s)_p\otimes
Y_{s,q}$ by the following formula: $$d(x\otimes y)=d(x)\otimes
y+(-1)^{n-s+p}x\otimes d(y).$$

A left differential module $(X,d,\mu)$ over a colored differential
algebra $(A,d,\pi)$ or, more briefly, a differential $A$-module is
a differential bigraded module $(X,d)$ endowed with a left action
$\mu:A\otimes X\to X$, which is a map of differential bigraded
modules of bidegree $(0,0)$ and satisfies the condition
$\mu(\pi\otimes 1)=\mu(1\otimes\mu)$.

Given a differential bigraded module $(X,d)$, consider the
corresponding colored differential algebra $(\hom(X;X),d,\pi)$. It
is easy to see that endowing $(X,d)$ with the structure of a left
differential module $(X,d,\mu)$ over a colored differential
algebra $(A,d,\pi)$ is equivalent to specifying a map
$\check{\mu}:(A,d,\pi)\to (\hom(X;X),d,\pi)$ of colored
differential algebras. Indeed, the maps $\mu$ and $\check{\mu}$
uniquely determine each other by the formula
$(\check{\mu}(a))(x)=\mu(a\otimes x)$, where $a\in A$, $x\in X$.

Now we note that for any quadratic colored algebra $(R,\pi)$ there
is the twisting cochain $\varphi_\vartheta^!:R^!_\bullet\to
R_{\bullet-1}$, which is defined by the following rule:
$$\varphi^!([x_1,\dots,x_k])=\left\{\begin{array}{ll}
x_1,&~\mbox{if}~k=1,\\ 0,&~\mbox{if}~k>0,\\
\end{array}\right.$$
where $(R^!,\nabla,\vartheta)$ is the colored coalgebra Koszul
dual to the quadratic colored algebra $(R,\pi)$.  It follows that,
as mentioned above, for any quadratic colored algebra $(R,\pi)$
the map of colored differential algebras
$\Omega(\varphi^!):(\Omega(R^!),d,\pi)\to (R,d=0,\pi)$ is defined,
which makes it possible to regard any differential $R$-module as a
differential $\Omega(R^!)$-module.

To describe the homotopy properties of differential
$\Omega(R^!)$-modules we first consider the colored variant of a
homotopy technique of differential Lie modules over graded
coalgebras \cite{Smir1}.

Let $(C,\nabla)$ be any connected colored graded coalgebra with
counit and coaugmentation.

A differential Lie module $(X,d,\psi)$ over a colored coalgebra
$(C,\nabla)$ or, briefly, a differential Lie $C$-module is an
arbitrary differential bigraded module $(X,d)$ equipped with the
map of bigraded modules $\psi:(C\otimes X)_{*,\bullet}\to
X_{*,\bullet-1}$ of bidegree $(0,-1)$ that satisfies the following
conditions:

1). $\psi(\nu\otimes 1)=d$, where $\nu:K_I\to C$ is the
coaugmentation of a colored graded coalgebra $(C,\nabla)$.

2). $d(\psi)+\psi\cup\psi=0$, where $\psi:(C\otimes X)_{*,\bu}\to
X_{*,\bu-1}$ is considered as the element $\psi\in\Hom(C\otimes
X;X)_{0,-1}$ of the differential bigraded module $(\Hom(C\otimes
X;X),d)$ and $\psi\cup\psi:(C\otimes X)_{*,\bu}\to X_{*,\bu-2}$ is
defined by the formula
$\psi\cup\psi=\psi(1\otimes\psi)(\nabla\otimes 1)$.

The same way as it was done for non-colored coalgebras in
\cite{Smir1}, it is easy to verify that the consideration of the
structure of a differential Lie module $(X,d,\psi)$ over a colored
graded coalgebra $(C,\nabla)$ on a differential bigraded module
$(X,d)$ is equivalent to the consideration of a perturbation
$t:(C\otimes X)_{*,\bu}\to (C\otimes X)_{*,\bu-1}$ of the
differential $1\otimes d:(C\otimes X)_{*,\bu}\to (C\otimes
X)_{*,\bu-1}$ on the differential module $(C\otimes X,1\otimes
d)$, which is a derivation of the free $C$-comodule $C\otimes X$.

A morphism $f:(X,d,\psi)\to (Y,d,\psi)$ of differential Lie
modules over a colored graded coalgebra $(C,\nabla)$ is a map of
bigraded modules $f:(C\otimes X)_{*,\bu}\to Y_{*,\bu}$ of bidegree
$(0,0)$ that satisfies the condition $d(f)+\psi\cup
f+f\cup\psi=0$, where $f$ is considered as the element
$f\in\Hom(C\otimes X;Y)_{0,0}$ of the differential bigraded module
$(\Hom(C\otimes X;X),d)$ and $\psi\cup f:(C\otimes X)_{*,\bu}\to
Y_{*,\bu-1}$, $f\cup\psi:(C\otimes X)_{*,\bu}\to Y_{*,\bu-1}$
respectively are defined by the formulae $$\psi\cup
f=\psi(1\otimes f)(\nabla\otimes 1),\quad
f\cup\psi=f(1\otimes\psi)(\nabla\otimes 1).$$

The composition $gf:(X,d,\psi)\to (X'',d,\psi)$ of morphisms
$f:(X,d,\psi)\to (X',d,\psi)$ and $g:(X',d,\psi)\to (X'',d,\psi)$
of differential Lie modules over a colored graded coalgebra
$(C,\nabla)$ is defined as the map $gf:(C\otimes X)_{*,\bu}\to
X''_{*,\bu}$ that is given by the formula $gf=g\cup f= g(1\otimes
f)(1\otimes\nabla)$.  The identity morphism $1_X:(X,d,\psi)\to
(X,d,\psi)$ for every differential Lie $C$-module $(X,d,\psi)$ is
defined as the map $1_X=\varepsilon\otimes 1 :(C\otimes
X)_{*,\bu}\to (K\otimes X)_{*,\bu}=X_{*,\bu}$, where
$\varepsilon:C\to K_I$ is counit of a colored graded coalgebra
$(C,\nabla)$.

A homotopy $h:(X,d,\psi)\to (Y,d,\psi)$ between morphisms
$f,g:(X,d,\psi)\to (Y,d,\psi)$ of differential Lie modules over a
colored graded algebra $(C,\nabla)$ is any map $h:(C\otimes
X)_{*,\bu}\to Y_{*,\bu+1}$ that satisfies the condition
$d(h)+\psi\cup h+h\cup\psi=f-g$, where $h$ is considered as the
element $h\in\Hom(C\otimes X;Y)_{0,1}$ of the differential
bigraded module $(\Hom(C\otimes X;Y),d)$ and maps $\psi\cup
h:(C\otimes X)_{*,\bu}\to Y_{*,\bu}$, $h\cup\psi:(C\otimes
X)_{*,\bu}\to Y_{*,\bu}$ respectively are defined by the formulae
$$\psi\cup h=\psi(1\otimes h)(\nabla\otimes
1),~~h\cup\psi=h(1\otimes\psi)(\nabla\otimes 1).$$

Let $\eta:(X,d,\psi)\,^{_{\longrightarrow}}_{^{\longleftarrow}}\,
(Y,d,\psi):\xi$ be any morphisms of differential Lie $C$-modules
such that $\eta\xi=1_Y$, and let $h:(X,d,\psi)\to (X,d,\psi)$ be a
homotopy between the morphisms $\xi\eta$ and $1_X$ of differential
Lie $C$-modules which satisfies the conditions $\eta h=0$, $\xi
h=0$, $h h=0$. Every triple
$(\eta:(X,d,\psi)\,^{_{\longrightarrow}}_{^{\longleftarrow}}\,
(Y,d,\psi):\xi\,,\,h)$ of the above form is said to be a SDR-data
of differential Lie $C$-modules.

The following assertion, which is proved in the same way as it is
done in \cite{Smir1} for non-colored coalgebras, establishes the
homotopy invariance of the structure of a differential Lie module
over any colored graded coalgebra under homotopy equivalences of
the type of SDR-data of differential bigraded modules.

{\bf Theorem 1.1}. Let $(X,d,\psi)$ be a differential Lie module
over a connected graded coalgebra $(C,\nabla)$ with counit
$\varepsilon:C\to K_I$ and coaugmentation $\nu:K_I\to C$, and let
$(\eta:(X,d)\,^{_{\longrightarrow}}_{^{\longleftarrow}}\,
(Y,d):\xi\,,\,h)$ be any SDR-data of differential bigraded
modules. Then the differential bigraded module $(Y,d)$ admits the
structure of a differential Lie $C$-module
$(Y,d,\overline{\psi})$, which is defined by the formula
$$\overline{\psi}=\varepsilon\otimes d+\eta t(1\otimes\xi)+\eta
t(1\otimes h)(1\otimes t)(\nabla\otimes
1)(1\otimes\xi)+\dots\eqno(1.1)$$ Moreover, the formulas
$$\overline{\xi}=\varepsilon\otimes\xi+(\varepsilon\otimes
h)(1\otimes t)(\nabla\otimes 1)(1\otimes\xi)\,+$$
$$+\,(\varepsilon\otimes h)(1\otimes t)(\nabla\otimes 1)(1\otimes
h)(1\otimes t)(\nabla\otimes 1)(1\otimes\xi)+\dots,\eqno(1.2)$$
$$\overline{\eta}=\varepsilon\otimes\eta+(\varepsilon\otimes\eta)(1\otimes
t)(\nabla\otimes 1)(1\otimes h)\,+$$
$$+\,(\varepsilon\otimes\eta)(1\otimes t)(\nabla\otimes
1)(1\otimes h)(1\otimes t)(\nabla\otimes 1)(1\otimes
h)+\dots,\eqno(1.3)$$ $$\overline{h}=\varepsilon\otimes
h+(\varepsilon\otimes h)(1\otimes t)(\nabla\otimes 1)(1\otimes
h)\,+$$ $$+\,(\varepsilon\otimes h)(1\otimes t)(\nabla\otimes
1)(1\otimes h)(1\otimes t)(\nabla\otimes 1)(1\otimes
h)+\dots\,,\eqno(1.4)$$ where $t=\psi-(\varepsilon\otimes d)$,
define the SDR-data
$(\overline{\eta}:(X,d,\psi)\,^{_{\longrightarrow}}_{^{\longleftarrow}}\,
(Y,d,\overline{\psi}):\overline{\xi}\,,\,\overline{h})$ of
differential Lie $C$-modules.~~~$\blacksquare$

Further the co-$B$-construction $(\Omega(R^!),d,\pi)$, where
$(R^!,\nabla)$ is the colored coalgebra Koszul dual to a quadratic
colored algebra $(R,\pi)$, we will be denoted by
$(R_\infty,d,\pi)$.

Let us consider homotopy properties of differential
$R_\infty$-modules. It is easy to see that the introduction of the
structure of a differential $R_\infty$-module on a differential
bigraded module $(X,d)$ is equivalent to the introduction of the
structure of a differential Lie $R^!$-module on $(X,d)$. Indeed,
every differential $R_\infty$-module $(X,d,\mu)$ defines the
differential Lie $R^!$-module $(X,d,\psi_\mu)$, where
$\psi_\mu=\mu(\varphi^\Omega\otimes 1):(R^!\otimes X)_{*,\bu}\to
X_{*,\bu-1}$. Conversely, the structure map $\psi:(R^!\otimes
X)_{*,\bu}\to X_{*,\bu-1}$ of a differential Lie $R^!$-module
$(X,d,\psi)$ determines the twisting cochain
$\check{\psi}:R^!_\bu\to\hom(X;X)_{\bu-1}$. It twisting cochain
induces the map
$\check{\mu}_\psi=\Omega(\check{\psi}):\Omega(R^!)=R_\infty\to
\hom(X;X)$ of colored differential algebras, which as said above
defines the structure map $\mu_\psi:R_\infty\otimes X\to X$ of the
differential $R_\infty$-module $(X,d,\mu_\psi)$.

{\bf Definition 1.3}.  By an $R_\infty$-map $f:(X,d,\mu)\to
(Y,d,\mu)$ of differential $R_\infty$-modules for a given any
quadratic colored algebra $(R,\pi)$ further we mean a morphism
$f:(X,d,\psi_\mu)\to (Y,d,\psi_\mu)$ of the corresponding
differential Lie $R^!$-modules.

Thus, an $R_\infty$-map $f:(X,d,\mu)\to (Y,d,\mu)$ of differential
$R_\infty$-modules is a map $f:(R^!\otimes X)_{*,\bu}\to
Y_{*,\bu}$ that satisfies the condition $d(f)+\psi_\mu\cup
f+f\cup\psi_\mu=0$.

{\bf Definition 1.4}. By an $R_\infty$-homotopy $h:(X,d,\mu)\to
(Y,d,\mu)$ between given $R_\infty$-maps $f,g:(X,d,\mu)\to
(Y,d,\mu)$ of differential $R_\infty$-modules for any quadratic
colored algebra $(R,\pi)$ we mean a homotopy $h:(X,d,\psi_\mu)\to
(Y,d,\psi_\mu)$ between morphisms $f,g:(X,d,\psi_\mu)\to
(Y,d,\psi_\mu)$ of the corresponding differential Lie
$R^!$-modules.

Thus, an $R_\infty$-homotopy $h:(X,d,\mu)\to (Y,d,\mu)$ between
given $R_\infty$-maps $f,g:(X,d,\mu)\to (Y,d,\mu)$ of differential
$R_\infty$-modules is a map $h:(R^!\otimes X)_{*,\bu}\to
Y_{*,\bu+1}$ that satisfies the condition $d(f)+\psi_\mu\cup
h+h\cup\psi_\mu=f-g$.

Suppose that
$\eta:(X,d,\mu)\,^{_{\longrightarrow}}_{^{\longleftarrow}}\,
(Y,d,\mu):\xi$ are $R_\infty$-maps of differential
$R_\infty$-modules such that $\eta\xi=1_Y$ and $h:(X,d,\mu)\to
(X,d,\mu)$ is an $R_\infty$-homotopy between the $R_\infty$-maps
$\xi\eta$ and $1_X$ of differential $R_\infty$-modules which
satisfies the conditions $\eta h=0$, $\xi h=0$, $hh=0$. Every
triple
$(\eta:(X,d,\mu)\,^{_{\longrightarrow}}_{^{\longleftarrow}}\,
(Y,d,\mu):\xi\,,\,h)$ of the above form is said to be a
$R_\infty$-SDR-data of differential $R_\infty$-modules.

The following theorem, which follows from Theorem 1.1, asserts the
$R_\infty$-homotopy invariance of the structure of a differential
$R_\infty$-module under homotopy equivalences of the type of
SDR-data of differential bigraded modules.

{\bf Theorem 1.2.}. Let an arbitrary differential
$R_\infty$-module $(X,d,\mu)$ and any SDR-data
$(\eta:(X,d)\,^{_{\longrightarrow}}_{^{\longleftarrow}}\,
(Y,d):\xi\,,\,h)$ of differential bigraded modules are given. Then
$(Y,d)$ can be equipped with the structure of a differential
$R_\infty$-module $(Y,d,\mu)$ given by the formula (1.1).
Moreover, there is an $R_\infty$-SDR-data
$(\overline{\eta}:(X,d,\mu)\,^{_{\longrightarrow}}_{^{\longleftarrow}}\,
(Y,d,\mu):\overline{\xi}\,,\,\overline{h})$ of differential
$R_\infty$-modules, which is defined by the formulas $(1.2)-(1.4)$
and satisfies the initial conditions $\overline{\eta}(\nu\otimes
1)=\eta$, $\overline{\xi}(\nu\otimes 1)=\xi$,
$\overline{h}(\nu\otimes 1)=h$, where $\nu:K_I\to R^!$ is
coaugmentation of the colored coalgebra $(R^!,\nabla)$ Koszul dual
to a quadratic colored algebra $(R,\pi)$.~~~$\blacksquare$
\vspace{0.5cm}

\centerline{\bf \S\,2. Differential modules with
$\infty$-simplicial faces}
\vspace{0.5cm}

Let $(F,\pi)$ be the colored graded $K$-algebra whose generators
are the elements $\p_i^n\in F(n-1,n)_0$, $n-1\in I$,
$i\in\mathbb{Z}$, $0\leqslant i\leqslant n$, connected by the
simplicial commutation relations
$$\p_i^{n-1}\p_j^n=\p_{j-1}^{n-1}\p_i^n,\quad i<j,~n-1\in I,$$
where $\p_i^{n-1}\p_j^n=\pi(\p_i^{n-1}\otimes \p_j^n)$. Further
the colored graded algebra $(F,\pi)$ will be called the colored
algebra of simplicial faces.

It is easy to see that the colored algebra of simplicial faces
$(F,\pi)$ is a quadratic colored algebra. Indeed, let $M$ be a
colored graded module defined by the conditions $M(s,t)_m=0$ for
$s,t\in I$, $m>0$, $M(s,t)_0=0$ for $(s,t)\not=(n-1,n)$, $n-1\in
I$, and $M(n-1,n)_0$ is a free $K$-module with the generators
$\p_i^n$, where $n-1\in I$, $0\leqslant i\leqslant n$. In the
colored graded module $M\otimes M$ we consider the submodule $Q$
that is defined by the conditions $Q(s,t)_m=0$ for $s,t\in I$,
$m>0$, $Q(s,t)_0=0$ for $(s,t)\not=(n-2,n)$, $n-2\in I$, and
$Q(n-2,n)_0$ is a free $K$-module with the generators
$\p_i^{n-1}\otimes\p_j^n-\p_{j-1}^{n-1}\otimes\p_i^n$, where
$n-2\in I$, $0\leqslant i<j\leqslant n$. It is clear that the
colored algebra of simplicial faces $(F,\pi)$ is isomorphic to the
quotient algebra $T(M)/(Q)$, and hence $(F,\pi)$ is a quadratic
colored algebra.

{\bf Definition 2.1}. By a differential module with simplicial
faces we mean an arbitrary differential module $(X,d,\mu)$ over
the colored algebra of simplicial faces $(F,\pi)$. By a map
$f:(X,d,\mu)\to (Y,d,\mu)$ of differential modules with simplicial
faces we mean any map $f:(X,d)\to (Y,d)$ of differential modules
that satisfies the condition $f\mu=\mu(1_F\otimes f)$. By a
homotopy between maps $f,g:(X,d,\mu)\to (Y,d,\mu)$ of differential
modules with simplicial faces we mean any homotopy $h:X_{*,\bu}\to
Y_{*,\bu+1}$ between maps $f,g:(X,d)\to (Y,d)$ of differential
bigraded modules that the condition $h\mu=\mu(1_F\otimes h)$
holds.

It is clear that the consideration of a differential module with
simplicial faces $(X,d,\mu)$ is equivalent to the consideration of
a differential module $(X,d)$ equipped with a family of the maps
$\check{\mu}(\p^n_i)=\p_i:X_{n,\,\bu}\to X_{n-1,\,\bu}$, where
$0\leqslant i\leqslant n$, $n\geqslant 0$, which are maps of
differential modules, i.e. satisfy the condition $d\p_i+\p_id=0$,
and also satisfy the simplicial commutation relations
$$\p_i\p_j=\p_{j-1}\p_i:X_{n,\bu}\to X_{n-2,\bu} ,\quad 0\leqslant
i<j\leqslant n.\eqno(2.1)$$ The above maps $\p_i$ of differential
modules are referred to as simplicial operators of faces or, more
briefly, simplicial faces of the differential bigraded module
$(X,d)$.

It is easy to see that the consideration of a map $f:(X,d,\p_i)\to
(Y,d,\p_i)$ of differential modules with simplicial faces is
equivalent to the consideration of a map $f:(X,d)\to (Y,d)$ of
differential modules that satisfies the condition
$$\p_if=f\p_i:X_{n,\bu}\to Y_{n-1,\bu},\quad 0\leqslant i\leqslant
n.\eqno(2.2)$$

The consideration of a homotopy between maps $f,g:(X,d,\p_i)\to
(Y,d,\p_i)$ of differential modules with simplicial faces is
equivalent to the consideration of a homotopy $h:X_{*,\bu}\to
Y_{*,\bu+1}$ between maps $f,g:(X,d)\to (Y,d)$ of differential
bigraded modules that satisfies the condition
$$\p_ih+h\p_i=0:X_{n,\bu}\to Y_{n-1,\bu+1},\quad 0\leqslant
i\leqslant n.\eqno(2.3)$$

Note that a difference between the concept of a differential
module with simplicial faces and the notion of a simplicial
differential module, i.e. of a simplicial object in the category
of differential modules \cite{May}, lies in the fact that for
differential modules with simplicial faces not assumed the
existence of simplicial degeneracy operators.

We now carry some necessary considerations for a introduction of
the concept of a differential module with $\infty$-simplicial
faces.

Let there be given arbitrary nonnegative integer $n\geqslant 0$.
Any collection of nonnegative integers $(i_1,\dots,i_k)$, where
$0\leqslant i_1<\dots<i_k\leqslant n$, further will be called a
ordered collection. Let $\Sigma_k$ be the symmetric group of
permutations of $k$ symbols. For any permutation
$\sigma\in\Sigma_k$ and any ordered collection $(i_1,\dots,i_k)$
we consider the collection $(\sigma(i_1),\dots,\sigma(i_k))$,
where $\sigma$ acts on the collection $(i_1,\dots,i_k)$ in a
standard way, i.e., permutes the numbers in this collection. For
this collection $(\sigma(i_1),\dots,\sigma(i_k))$ we define the
collection $(\widehat{\sigma(i_1)},\dots,\widehat{\sigma(i_k)})$
by setting $$\widehat{\sigma
(i_s)}=\sigma(i_s)-\alpha(\sigma(i_s)),\quad 1\leqslant s\leqslant
k,$$ where $\alpha(\sigma(i_s))$ is the number of those elements
of $(\sigma(i_1),\dots,\sigma(i_s),\dots\sigma(i_k))$ on the right
of $\sigma(i_s)$ which are smaller than $\sigma(i_s)$.

Consider the colored coalgebra $(F^!,\nabla)$ Koszul dual to the
colored algebra of simplicial faces $(F,\pi)$. It is easy to see
that generators of the free $K$-module $(F^!)^{(k)}(n-k,n)_k$,
$k\geqslant 1$, $n-k\in I$, are the elements
$$[\p_{i_1}^{n-k+1}]\hat{\wedge}\dots\hat{\wedge}\,
[\p_{i_k}^n]=\sum_{\sigma\in \Sigma_k}(-1)^{{\rm
sign}(\sigma)+1}[\p_{\widehat{\sigma(i_1)}}^{n-k+1},\dots,\p_{\widehat{\sigma(i_k)}}^n]
,~~0\leqslant i_1<\dots<i_k\leqslant n,$$ and that the equality
$(F^!)(k)(s,t)_m=0$ for $(s,t)\not=(n-k,n)$, $m\not=k$, $n-k\in
I$, holds.

A direct calculation shows that a comultiplication $\nabla$ of the
colored graded coalgebra $F^!$ on the generators
$[\p_{i_1}^{n-k+1}]\hat{\wedge}\dots\hat{\wedge}\,[\p_{i_k}^n]$ of
the module $(F^!)^{(k)}(n-k,n)_k$, $k\geqslant 1$, $n-k\in I$, is
given by the following formula:
$$\nabla([\p_{i_1}^{n-k+1}]\hat{\wedge}\dots\hat{\wedge}\,[\p_{i_k}^n])=$$
$$=1_{n-k}\otimes
([\p_{i_1}^{n-k+1}]\hat{\wedge}\dots\hat{\wedge}\,[\p_{i_k}^n])+
([\p_{i_1}^{n-k+1}]\hat{\wedge}\dots\hat{\wedge}\,[\p_{i_k}^n])\otimes
1_n+$$ $$+\sum_{\sigma\in\Sigma_k}\sum_{I_{\sigma}} (-1)^{{\rm
sign}(\sigma)+1}
([\p_{\widehat{\sigma(i_1)}}^{n-k+1}]\hat{\wedge}\dots\hat{\wedge}\,[\p_{\widehat{\sigma(i_m)}}^{n-k+m}])\otimes
([\p_{(\widehat{\sigma(i_{m+1})}}^{n-k+m+1}]\hat{\wedge}\dots\hat{\wedge}\,[\p_{\widehat{\sigma(i_k)}}^n]),$$
where  $I_\sigma$ is the set all partitions of the collection
$(\widehat{\sigma(i_1)},\dots,\widehat{\sigma(i_k)})$ into two
ordered collections
$(\widehat{\sigma(i_1)},\dots,\widehat{\sigma(i_m)})$ and
$(\widehat{\sigma(i_{m+1})},\dots,\widehat{\sigma(i_k)})$,
$1\leqslant m\leqslant k-1$.

{\bf Definition 2.2}. The colored differential algebra
$(\Omega(F^!),d,\pi)$ will be called the colored algebra of
$\infty$-simplicial faces and denoted by $(F_\infty,d,\pi)$.

It is easy to see that generators of colored algebra $(\F,\pi)$
are the elements
$\p^n_{(i_1,\dots,i_k)}=[[\p_{i_1}^{n-k+1}]\hat{\wedge}\dots
\hat{\wedge}\,[\p_{i_k}^n]]\in\F(n-k,n)_{k-1}$, $0\leqslant
i_1<\dots<i_k\leqslant n$, $n-k\in I$, that connected by the
relations
$$d(\p^n_{(i_1,\dots,i_k)})=\sum_{\sigma\in\Sigma_k}\sum_{I_{\sigma}}
(-1)^{{\rm sign}(\sigma)+1}
\p^{n-(k-m)}_{(\widehat{\sigma(i_1)},\dots,\widehat{\sigma(i_m)})}
\p^n_{(\widehat{\sigma(i_{m+1})},\dots,\widehat{\sigma(i_k)})},\eqno(2.4)$$
where $\p_{(\widehat{\sigma(i_1)},\dots,\widehat{\sigma(i_m)})}
\p_{(\widehat{\sigma(i_{m+1})},\dots,\widehat{\sigma(i_k)})}=
\pi(\p_{(\widehat{\sigma(i_1)},\dots,\widehat{\sigma(i_m)})}
\otimes\p_{(\widehat{\sigma(i_{m+1})},\dots,\widehat{\sigma(i_k)})})$,
and  $I_{\sigma}$ is just like that in the above formula for the
comultiplication of the colored graded coalgebra $(F^!,\nabla)$.

The following statement helps to simplify the procedure for
enumerating all partitions from the set $I_\sigma$, by which goes
a summation in the formula $(2.4)$.

{\bf Proposition 2.1}. Let any ordered collection
$(i_1,\dots,i_k)$ and any permutation $\sigma\in\Sigma_k$ are
given. The partition
$(\sigma(i_1),\dots,\sigma(i_m)|\,\sigma(i_{m+1}),\dots,\sigma(i_k))$
is a partition into two ordered collections if and only if the
partition
$(\widehat{\sigma(i_1)},\dots,\widehat{\sigma(i_m)}\,|\,\widehat{\sigma(i_{m+1})},
\dots,\widehat{\sigma(i_k)})$ also is a partition into two ordered
collections.

{\bf Proof}. Suppose that
$(\sigma(i_1),\dots,\sigma(i_m)\,|\,\sigma(i_{m+1}),\dots,\sigma(i_k))$
is a partition into two ordered collections. Let us show that
$(\widehat{\sigma(i_1)},\dots,\widehat{\sigma(i_m)}\,|\,\widehat{\sigma(i_{m+1})},
\dots,\widehat{\sigma(i_k)})$ also is a partition into two ordered
collections. It is easy to see that
$\widehat{\sigma(i_t)}=\sigma(i_t)$ for $m+1\leqslant t\leqslant
k$ and hence we have
$\widehat{\sigma(i_{m+1})}<\dots<\widehat{\sigma(i_k)}$. Denote by
$\beta(\sigma(i_s))$, where $2\leqslant s\leqslant m$, the number
of those elements of
$(\sigma(i_1),\dots,\sigma(i_s),\dots,\sigma(i_k))$ on the right
of $\sigma(i_s)$ which are smaller than $\sigma(i_s)$ and larger
than $\sigma(i_{s-1})$. It is clear that
$\beta(\sigma(i_s))=\alpha(\sigma(i_s))-\alpha(\sigma(i_{s-1}))$.
Moreover, since $\sigma(i_{s-1})<\sigma(i_s)$, it follows that
$\beta(\sigma(i_s))\leqslant\sigma(i_s)-\sigma(i_{s-1})-1$. From
this we obtain
$$\widehat{\sigma(i_{s-1})}=\sigma(i_{s-1})-\alpha(\sigma(i_{s-1}))
<\sigma(i_s)-\beta(\sigma(i_s))-\alpha(\sigma(i_{s-1}))=\widehat{\sigma(i_s)}$$
and hence we have
$\widehat{\sigma(i_1)}<\dots<\widehat{\sigma(i_m)}$. We now show
the converse. Suppose that
$(\widehat{\sigma(i_1)},\dots,\widehat{\sigma(i_m)}\,|\,\widehat{\sigma(i_{m+1})},
\dots,\widehat{\sigma(i_k)})$ is a partition into two ordered
collections. Since
$\widehat{\sigma(i_{s-1})}<\widehat{\sigma(i_s)}$, $2\leqslant
s\leqslant m$, it follows that
$\sigma(i_{s-1})-\sigma(i_s)<\alpha(\sigma(i_{s-1}))-\alpha(\sigma(i_s))$.
Now assume the opposite, i.e. assume that for some $2\leqslant
s\leqslant m$ the condition $\sigma(i_{s-1})>\sigma(i_s)$ is true.
For the above $s$ we denote by $\gamma(\sigma(i_{s-1}))$ the
number of those elements of
$(\sigma(i_1),\dots,\sigma(i_s),\dots,\sigma(i_k))$ on the right
of $\sigma(i_{s-1})$ which are smaller than $\sigma(i_{s-1})$ and
larger than $\sigma(i_s)$ or equal to $\sigma(i_s)$. It is easy to
see that
$\alpha(\sigma(i_{s-1}))=\alpha(\sigma(i_s))+\gamma(\sigma(i_{s-1}))$
and $\gamma(\sigma(i_{s-1}))\leqslant\sigma(i_{s-1})-\sigma(i_s)$.
From this we obtain the condition
$$\alpha(\sigma(i_{s-1}))-\alpha(\sigma(i_s))=\gamma(\sigma(i_{s-1}))
\leqslant\sigma(i_{s-1})-\sigma(i_s),$$ which contradicts to the
condition
$\sigma(i_{s-1})-\sigma(i_s)<\alpha(\sigma(i_{s-1}))-\alpha(\sigma(i_s))$.
Thus, we have $\sigma(i_1)<\dots<\sigma(i_m)$. Carrying out
similar arguments in the case
$\widehat{\sigma(i_{s-1})}<\widehat{\sigma(i_s)}$, $m+2\leqslant
s\leqslant k$, we obtain
$\sigma(i_{m+1})<\dots<\sigma(i_k)$.~~~$\blacksquare$

{\bf Definition 2.3}. By a differential module with
$\infty$-simplicial faces we mean an arbitrary differential module
$(X,d,\mu)$ over the colored algebra of $\infty$-simplicial faces
$(F_\infty,d,\pi)$. By a morphism $f:(X,d,\mu)\to (Y,d,\mu)$ of
differential modules with $\infty$-simplicial faces we mean any
$F_\infty$-map $f:(X,d,\mu)\to (Y,d,\mu)$ of differential
$F_\infty$-modules. By a homotopy between morphisms
$f,g:(X,d,\mu)\to (Y,d,\mu)$ of differential modules with
$\infty$-simplicial faces we mean any $F_\infty$-homotopy
$h:(X,d,\mu)\to (Y,d,\mu)$ between $F_\infty$-maps $f$ and $g$ of
differential $F_\infty$-modules. By a SDR-data of differential
modules with $\infty$-simplicial faces we mean any $\F$-SDR-Data
of differential $F_\infty$-modules.

It is easy to see that the consideration of a differential module
with $\infty$-simplicial faces $(X,d,\mu)$ is equivalent to the
consideration of a differential bigraded module $(X,d)$ equipped
with a family of module maps
$$\widetilde{\p}=\{\p_{(i_1,\dots,i_k)}=\check{\mu}([[\p_{i_1}^{n-k+1}]\hat{\wedge}\dots
\hat{\wedge}\,[\p_{i_k}^n]]):X_{n,\bu}\to X_{n-k,\bu+k-1}\},\quad
n\geqslant 0,$$ where $0\leqslant i_1<\dots<i_k\leqslant n$, which
satisfy the following relations:
$$d(\p_{(i_1,\dots,i_k)})=\sum_{\sigma\in\Sigma_k}\sum_{I_{\sigma}}
(-1)^{{\rm sign}(\sigma)+1}
\p_{(\widehat{\sigma(i_1)},\dots,\widehat{\sigma(i_m)})}
\p_{(\widehat{\sigma(i_{m+1})},\dots,\widehat{\sigma(i_k)})}.\eqno(2.5)$$
The above module maps $\p_{(i_1,\dots,i_k)}$ will be called
$\infty$-simplicial faces of a differential bigraded module
$(X,d)$. Further, each differential module with
$\infty$-simplicial faces $(X,d,\mu)$ will be identified with the
corresponding triple $(X,d,\widetilde{\p})$.

For $k=1$, the relations $(2.5)$ take the form $d(\p_{(i)})=0$.
This indicates that all $\infty$-simplicial faces $\p_{(i)}$ are
maps of differential modules. For $k=2$, relations $(2.5)$ take
the form $$d(\p_{(i,j)})=\p_{(j-1)}\p_{(i)}-\p_{(i)}\p_{(j)},\quad
i<j.$$ This means that the $\infty$-simplicial face $\p_{(i,j)}$
is a homotopy between the maps of differential modules
$\p_{(j-1)}\p_{(i)}$ and $\p_{(i)}\p_{(j)}$. Thus,
$\infty$-simplicial faces $\p_{(i)}$ satisfy the simplicial
commutation relations $(2.1)$ up to homotopy. For $k=3$, by using
Proposition 2.1 we easily obtain that the relations $(2.5)$ take
the following form:
$$d(\p_{(i_1,i_2,i_3)})=-\p_{(i_1)}\p_{(i_2,i_3)}-\p_{(i_1,i_2)}\p_{(i_3)}-
\p_{(i_3-2)}\p_{(i_1,i_2)}\,-$$
$$-\,\p_{(i_2-1,i_3-1)}\p_{(i_1)}+\p_{(i_2-1)}\p_{(i_1,i_3)}+\p_{(i_1,i_3-1)}\p_{(i_2)},
\qquad i_1<i_2<i_3.$$

The consideration of a morphism $f:(X,d,\mu)\to (Y,d,\mu)$ of
differential modules with $\infty$-simplicial faces is equivalent
to the consideration of a family of module maps
$$\widetilde{f}=\{f_{(\,\,)}=\check{f}(1_n):X_{n,\bu}\to
Y_{n,\bu},$$
$$f_{(i_1,\dots,i_k)}=\check{f}([\p_{i_1}]\hat{\wedge}\dots\hat{\wedge}\,[\p_{i_k}]):X_{n,\bu}\to
Y_{n-k,\bu+k}\},\quad n\geqslant 0,$$ where $1_n=1\in
K=(F^!)^{(0)}(n,n)_0$, $\check{f}(a)(x)=f(a\otimes x)$, $a\in
F^!$, $0\leqslant i_1<\dots<i_k\leqslant n$, which satisfy the
conditions $$d(f_{(\,\,)})=0,\quad
d(f_{(i_1,\dots,i_k)})=-\p_{(i_1,\dots,i_k)}f_{(\,\,)}+f_{(\,\,)}\p_{(i_1,\dots,i_k)}+\eqno(2.6)$$
$$=\sum_{\sigma\in\Sigma_k}\sum_{I_{\sigma}}(-1)^{{\rm
sign}(\sigma)+1}\p_{(\widehat{\sigma(i_1)},\dots,
\widehat{\sigma(i_m)})}f_{(\widehat{\sigma(i_{m+1})},\dots,
\widehat{\sigma(i_k)})}-f_{(\widehat{\sigma(i_1)},\dots,
\widehat{\sigma(i_m)})}\p_{(\widehat{\sigma(i_{m+1})},\dots,
\widehat{\sigma(i_k)})}.$$ Further, each morphism $f:(X,d,\mu)\to
(Y,d,\mu)$ of differential modules with $\infty$-simplicial faces
will be identified with the corresponding family of maps
$\widetilde{f}:(X,d,\widetilde{\p})\to (Y,d,\widetilde{\p})$.

For $k=1$, the relations $(2.6)$ take the form
$$d(f_{(i)})=f_{(\,\,)}\p_{(i)}-\p_{(i)}f_{(\,\,)},\qquad
i\geqslant 0.$$ This means that the map $f_{(i)}$ is a homotopy
between the maps of differential modules $f_{(\,\,)}\p_{(i)}$ and
$\p_{(i)}f_{(\,\,)}$. Thus, the map of differential modules
$f_{(\,\,)}$ satisfies the condition $(2.2)$ up to homotopy. For
$k=2$, the relations $(2.6)$ take the following form:
$$d(f_{(i,j)})=-\p_{(i,j)}f_{(\,\,)}+f_{(\,\,)}\p_{(i,j)}-\p_{(i)}f_{(j)}+
\p_{(j-1)}f_{(i)}+ f_{(i)}\p_{(j)}-f_{(j-1)}\p_{(i)},\quad i<j.$$

The consideration of a homotopy $h:(X,d,\mu)\to (Y,d,\mu)$ between
morphism $f,g:(X,d,\mu)\to (Y,d,\mu)$ of differential modules with
$\infty$-simplicial faces is equivalent to the consideration of a
family of module maps
$$\widetilde{h}=\{h_{(\,\,)}=\check{h}(1_n):X_{n,\bu}\to
Y_{n,\bu+1},$$
$$h_{(i_1,\dots,i_k)}=\check{h}([\p_{i_1}]\hat{\wedge}\dots\hat{\wedge}\,[\p_{i_k}]):X_{n,\bu}\to
Y_{n-k,\bu+k+1}\},\quad n\geqslant 0,$$ where $1_n=1\in
K=(F^!)^{(0)}(n,n)_0$, $\check{h}(a)(x)=h(a\otimes x)$, $a\in
F^!$, $0\leqslant i_1<\dots<i_k\leqslant n$, which satisfy the
conditions
$$d(h_{(\,\,)})=f_{(\,\,)}-g_{(\,\,)},~~d(h_{(i_1,\dots,i_k)})=\eqno(2.7)$$
$$=f_{(i_1,\dots,i_k)}-g_{(i_1,\dots,i_k)}-
\p_{(i_1,\dots,i_k)}h_{(\,\,)}-h_{(\,\,)}\p_{(i_1,\dots,i_k)}+$$
$$+\sum_{\sigma\in\Sigma_k}\sum_{I_{\sigma}}(-1)^{{\sign(\sigma)}+1}
\p_{(\widehat{\sigma(i_1)},\dots,\widehat{\sigma(i_m)})}
h_{(\widehat{\sigma(i_{m+1})},\dots,\widehat{\sigma(i_k)})}+
h_{(\widehat{\sigma(i_1)},\dots,\widehat{\sigma(i_m)})}
\p_{(\widehat{\sigma(i_{m+1})},\dots,\widehat{\sigma(i_k)})}.$$
Further, each homotopy $h:(X,d,\mu)\to (Y,d,\mu)$ between
morphisms $f,g:(X,d,\mu)\to (Y,d,\mu)$ of differential modules
with $\infty$-simplicial faces will be identified with the
corresponding family of maps
$\widetilde{h}:(X,d,\widetilde{\p})\to (Y,d,\widetilde{\p})$.

For $k=1$, the relations $(2.7)$ take the form
$$d(h_{(i)})=f_{(i)}-g_{(i)}-\p_{(i)}h_{(\,\,)}-h_{(\,\,)}\p_{(i)},\qquad
i\geqslant 0.$$

These equalities admit a good interpretation in the following
special case. Suppose that $f,g:(X,d,\p_i)\to (Y,d,\p_i)$ are any
maps of differential modules with simplicial faces, which we
regard as morphisms $\widetilde{f},\widetilde{g}:X\to Y$ of
differential modules with $\infty$-simplicial faces, and let
$\widetilde{h}:X\to Y$ be a homotopy between $\widetilde{f}$ and
$\widetilde{g}$. In the case under consideration, we have
$f_{(i)}=g_{(i)}=0$; therefore, the relations
$d(h_{(i)})=f_{(i)}-g_{(i)}-\p_{(i)}h_{(\,\,)}-h_{(\,\,)}\p_{(i)}$,
$i\geqslant 0$, can be written as the relations
$$d(h_{(i)})=0-(\p_{(i)}h_{(\,\,)}+h_{(\,\,)}\p_{(i)}),~~i\geqslant
0,$$ which means that the map $h_{(i)}$ is a homotopy between the
maps $0$ and $\p_{(i)}h_{(\,\,)}+h_{(\,\,)}\p_{(i)}$ of
differential modules. Thus, in considered case the homotopy
$h_{(\,\,)}$ between $f_{(\,\,)}$ and $g_{(\,\,)}$ satisfies the
condition $(2.3)$ up to homotopy. For $k=2$, the relations $(2.7)$
take the following form: $$d(h_{(i,j)})=f_{(i,j)}-g_{(i,j)}-
\p_{(i,j)}h_{(\,\,)}-h_{(\,\,)}\p_{(i,j)}-\p_{(i)}h_{(j)}\,+$$
$$+\,\p_{(j-1)}h_{(i)}- h_{(i)}\p_{(j)}+h_{(j-1)}\p_{(i)},\qquad
i<j.$$

Applying the theorem 1.2 to the colored algebra $(F,\pi)$ of
simplicial faces we obtain the assertion, which  establishes the
homotopy invariance of the structure of a differential module with
$\infty$-simplicial faces under homotopy equivalences of the type
of SDR-data of differential bigraded modules.

{\bf Theorem 2.1}. Let any differential module with
$\infty$-simplicial faces $(X,d,\widetilde{\p})$ and an arbitrary
SDR-data $(\eta:(X,d)\,^{_{\longrightarrow}}_{^{\longleftarrow}}\,
(Y,d):\xi\,,\,h)$ of differential bigraded modules are given. Then
$(Y,d)$ can be equipped with the structure of a differential
module with $\infty$-simplicial faces $(Y,d,\widetilde{\p})$ given
by the formula (1.1). Moreover, there is an $F_\infty$-SDR-data
$(\overline{\eta}:(X,d,\mu)\,^{_{\longrightarrow}}_{^{\longleftarrow}}\,
(Y,d,\mu):\overline{\xi}\,,\,\overline{h})$ of differential
modules with $\infty$-simplicial faces, which is defined by the
formulas $(1.2)-(1.4)$ and satisfies the initial conditions
$\eta_{(\,\,)}=\eta$, $\xi_{(\,\,)}=\xi$,
$h_{(\,\,)}=h$.~~~$\blacksquare$ \vspace{0.5cm}

\centerline{\bf \S\,3. The chain realization of differential
modules with} \centerline{\bf $\infty$-simplicial faces and the
$B$-construction for $\A$-algebras}
\vspace{0.5cm}

For the colored algebra of $\infty$-simplicial faces $(\F,d,\pi)$
the differential graded module $(\F(m,n)_\bu,d)$ we denoted by
$(\F[n]_{m,\bu},d)$. It is easy to see that for each $n\geqslant
0$ the differential bigraded module $(\F[n],d)$,
$\F[n]=\{\F[n]_{m,p}\}$, $m\geqslant 0$, $p\geqslant 0$,
$d:F[n]_{m,\bu}\to F[n]_{m,\bu-1}$, equipped with the structure of
a differential module with $\infty$-simplicial faces
$\widetilde{\p}=\{\p_{(i_1,\dots,i_k)}:\F[n]_{m,\bu}\to
\F[n]_{m-k,\bu+k-1}\}$, which for every element $a\in
\F[n]_{m,\bu}=\F(m,n)_\bu$ are defined by the following equality:
$$\p_{(i_1,\dots,i_k)}(a)=\pi(\p^m_{(i_1,\dots,i_k)}\otimes
a),\quad \p_{(i_1,\dots,i_k)}^m\in \F(m-k,m)_{k-1}.$$ It is clear
that for the differential module with $\infty$-simplicial faces
$(\F[n],d,\widetilde{\p})$ its corresponding differential module
$(\F[n],d,\mu)$ over the colored algebra of $\infty$-simplicial
faces $(\F,d,\pi)$ is the free bigraded module $(\F[n],\mu)$ over
the colored graded algebra $(\F,\pi)$ with one generator $1_n\in
\F[n]_{n,0}=\F(n,n)_0$, and the differential $d:F[n]_{m,\bu}\to
F[n]_{m,\bu-1}$ completely defined by the differential in
$(\F,d,\pi)$ and the equality $d(1_n)=0$, i.e. is defined for any
elements $\mu(a\otimes 1_n)\in\F[n]$ by the formula
$d(\mu(a\otimes 1_n))=\mu(d(a)\otimes 1_n)$.

Now suppose that an arbitrary differential module with
$\infty$-simplicial faces $(X,d,\widetilde{\p})$ is given.
Consider for any fixed element $x\in X_{n,q}$ the map
$\overline{x}:\F[n]_{*,\bu}\to X_{*,\bu+q}$ of bigraded modules,
which for every element $a\in \F[n]_{m,p}$ is defined by the
formula $$\overline{x}(a)=(-1)^{(p+m)q}\mu(a\otimes
x),\eqno(3.1)$$ where $(X,d,\mu)$ is the differential
$F_\infty$-module that corresponds to the differential module with
$\infty$-simplicial faces $(X,d,\widetilde{\p})$. By the
definition of a differential in the tensor product of a colored
differential module and the differential bigraded module we obtain
that the map $\overline{x}:\F[n]_{*,\bu}\to X_{*,\bu+q}$, which
considered as the element $\overline{x}\in\Hom(\F[n];X)_{0,q}$ of
the differential bigraded module $(\Hom(\F[n];X),d)$, satisfies
the following formula:
$$d(\overline{x})=(-1)^n\overline{d(x)}.\eqno(3.2)$$ This formula
implies that for an arbitrary cycle $x\in X_{n,q}$ of the
differential bigraded module $(X,d)$ its corresponding map
$\overline{x}$ is the map of differential bigraded modules
$\overline{x}:(\F[n]_{*,\bu},d)\to (X_{*,\bu+q},d)$ of bidegree
$(0,q)$, i.e. the map $\overline{x}$ satisfies the condition
$d\overline{x}=(-1)^q\overline{x}d$.

Let us consider for an arbitrary differential $\F$-module
$(X,d,\mu)$ and any elements $a\in\F(s,m)_t=\F[m]_{s,t}$ and $x\in
X_{m,p}$ the element $\mu(a\otimes x)\in X_{s,t+p}$ and the maps
$\overline{a}:\F[s]_{*,\bu}\to\F[m]_{*,\bu+t}$,
$\overline{x}:\F[m]_{*,\bu}\to X_{*,\bu+p}$,
$\overline{\mu(a\otimes x)}:\F[s]_{*,\bu}\to X_{*,\bu+t+p}$, which
correspond to the above elements. It is easy to see that $(3.1)$
implies the following true relation: $$\overline{\mu(a\otimes
x)}=(-1)^{tp}\overline{x}\,\overline{a}.\eqno(3.3)$$ In
particular, if as a differential $\F$-module $(X,d,\mu)$ is take
the differential $\F$-module $(\F[n],d,\mu)$, then for the above
element $\mu(a\otimes x)=\pi(a\otimes x)=ax$ the relation $(3.3)$
can be written as the relation
$$\overline{ax}=(-1)^{tp}\overline{x}\,\overline{a}.\eqno(3.3\,')$$

By using the formula $(3.1)$ we obtain that for any element $x\in
X_{n,q}$ its corresponding map $\overline{x}:\F[n]_{*,\bu}\to
X_{*,\bu+q}$ is related to $\infty$-simplicial faces
$\p_{(i_1,\dots,i_k)}$ on $\F[n]$ and on $X$ by the formula
$$\p_{(i_1,\dots,i_k)}\overline{x}=(-1)^q\overline{x}\,\p_{(i_1,\dots,i_k)}.\eqno(3.4)$$
This formula implies that for an arbitrary element $x\in X_{n,q}$
its corresponding map $\overline{x}:\F[n]_{*,\bu}\to X_{*,\bu+q}$
is completely defined by its value on the generator $1_n\in
\F[n]_{n,0}$, i.e. completely defined by the equality
$\overline{x}(1_n)=(-1)^{nq}x$.

Now, let us consider for each element $\p^n_{(i_1,\dots,i_k)}\in
\F[n]_{n-k,k-1}$ its corresponding map
$$\delta^{(i_1,\dots,i_k)}=(-1)^{n-k}\overline{\p^n_{(i_1,\dots,i_k)}}:\F[n-k]_{*,\bu}\to
\F[n]_{*,\bu+k-1}.$$ Further the maps $\delta^{(i_1,\dots,i_k)}$
we will be called $\infty$-cosimplicial cofaces of the family
$\F[*]=\{(\F[n],d,\widetilde{\p})\}_{n\geqslant 0}$ of
differential modules with $\infty$-simplicial faces. By using the
formulas $(2.4)$, $(3.2)$ and $(3.3\,')$ we obtain that
$\infty$-cosimplicial cofaces of the family $\F[*]$ related by
following relations:
$$d(\delta^{(i_1,\dots,i_k)})=\sum_{\sigma\in\Sigma_k}\sum_{I_{\sigma}}
(-1)^{{\rm
sign}(\sigma)+k(m-1)}\delta^{(\widehat{\sigma(i_{m+1})},\dots,\widehat{\sigma(i_k)})}
\delta^{(\widehat{\sigma(i_1)},\dots,\widehat{\sigma(i_m)})},\eqno(3.5)$$
where the set $I_{\sigma}$ is just that in the formula $(2.4)$.

Recall \cite{Lap1} that by a $\D$-module one means any bigraded
module $X=\{X_{n,m}\}$, $n,m\in\mathbb{Z}$, equipped with a family
$\{d^k:X_{*,\bu} \to
X_{*-k,\bu+k-1}\,\,|\,\,k\in\mathbb{Z},\,\,k\geqslant 0\}$ of
module maps, which for each integer $k\geqslant 0$ satisfy the
following relation: $$\sum\limits_{i+j=k}d^id^j=0.\eqno(3.6)$$

For every differential module with $\infty$-simplicial faces
$(X,d,\widetilde{\p})$, we define the family of maps
$\{d^k:X_{n,\bu}\to X_{n-k,\bu+k-1}\}$, $k,n\in\mathbb{Z}$,
$k\geqslant 0$, $n\geqslant 0$, by setting $$\left.
\begin{array}{c}d^0=d:X_{n,\bu}\to X_{n,\bu-1},\quad d^k=0:X_{n,\bu}\to
X_{n-k,\bu+k-1},\quad k>n,\\
\\
d^k=\sum\limits_{0\leqslant i_1<\dots<i_k\leqslant n
}(-1)^{i_1+\dots+i_k}\p_{(i_1,\dots,i_k)}:X_{n,\bu}\to
X_{n-k,\bu+k-1},\quad k\leqslant n\\
\end{array}\right\}\eqno(3.7)$$
Since for any permutation $\sigma\in\Sigma_k$ of integers
$i_1<\dots<i_k$ there is the true equality
$$\widehat{\sigma(i_1)}+\dots+\widehat{\sigma(i_k)}=
\sigma(i_1)+\dots+\sigma(i_k)-I(\sigma),$$ where
$I(\sigma)=\alpha(\sigma(i_1))+\dots+\alpha(\sigma(i_k))$ is a
number of the permutation $\sigma$, by using
${\sign}(\sigma)\equiv I(\sigma){\rm mod}(2)$ and $(2.5)$, it is
easy to check that, for maps from the family $(3.7)$, the
relations $(3.6)$ are true.

Now consider for each differential module with $\infty$-simplicial
faces $(X,d,\widetilde{\p})$ its corresponding $\D$-module
$(X,d^{\,k})$ and define the differential graded module
$(\overline{X},\p)$ by setting
$$\overline{X}_n=\bigoplus_{s+t=n}X_{s,t},\quad
\p=\sum_{k\geqslant 0}d^{\,k}:\overline{X}_\bu\to
\overline{X}_{\bu-1}.$$ It is easy to see that from the relations
$(3.6)$ follows the equality $\p\p=0$. In particular, for each
differential module with $\infty$-simplicial faces
$(\F[n],d,\widetilde{\p})$, $n\geqslant 0$, there is the
differential module $(\overline{\F[n]},\p)$.

The above considered map $\overline{x}:\F[n]_{*,\bu}\to
X_{*,\bu+q}$ that corresponds to an arbitrary element $x\in
X_{n,q}$ induces the map of graded modules
$\overline{x}:\overline{\F[n]}_\bu\to \overline{X}_{\bu+q}$ of
degree $q$, which, by the equality $(3.4)$, satisfies the relation
$\p(\overline{x})=d(\overline{x})$. In particular, since
$\p(\delta^{(i_1,\dots,i_k)})=d(\delta^{(i_1,\dots,i_k)})\not=0$,
the maps $\delta^{(i_1,\dots,i_k)}:\overline{\F[n-k]}_\bu\to
\overline{\F[n]}_{\bu+k-1}$, which are induced by
$\infty$-cosimplicial cofaces
$\delta^{(i_1,\dots,i_k)}:\F[n-k]_{*,\bu}\to \F[n]_{*,\bu+k-1}$,
$k\geqslant 2$, are not maps of differential modules.

Now suppose that an arbitrary differential module with
$\infty$-simplicial faces $(X,d,\widetilde{\p})$ is given. For
each fixed $n\geqslant 0$, let us consider the tensor product
$(\overline{\F[n]}_\bu\otimes X_{n,\bu},d)$ of differential
modules $(\overline{\F[n]}_\bu,\p)$ and $(X_{n,\bu},d)$.

{\bf Theorem 3.1}. The differential in the direct sum
$\bigoplus_{n\geqslant 0}(\overline{\F[n]}_\bu\otimes
X_{n,\bu},d)$ of differential modules induces a well defined
differential in the quotient module $$\bigoplus_{n\geqslant
0}(\overline{\F[n]}_\bu\otimes X_{n,\bu})/\!\sim$$ of the graded
module $\bigoplus_{n\geqslant 0}(\overline{\F[n]}_\bu\otimes
X_{n,\bu})$ by the equivalence relation $\sim\,$, which is
generated by the following relations:
$$\delta^{(i_1,\dots,i_k)}(a)\otimes x\sim
(-1)^{s(k-1)+(n+q)k}a\otimes\p_{(i_1,\dots,i_k)}(x),\eqno(3.8)$$
where $a\in\overline{\F[n-k]}_s$, $x\in X_{n,q}$ are an arbitrary
elements, and $(i_1,\dots,i_k)$ is any ordered collection.

{\bf Proof}. It suffices to show that from the condition $(3.8)$
follows the condition $d(\delta^{(i_1,\dots,i_k)}(a)\otimes x)\sim
d((-1)^{s(k-1)+(n+q)k}a\otimes\p_{(i_1,\dots,i_k)}(x))$. Since
$$\p\delta^{(i_1,\dots,i_k)}+(-1)^k\delta^{(i_1,\dots,i_k)}\p=
\p(\delta^{(i_1,\dots,i_k)})=d(\delta^{(i_1,\dots,i_k)}),$$ by
using the formula $(3.5)$ we obtain
$$d(\delta^{(i_1,\dots,i_k)}(a)\otimes
x)=\p(\delta^{(i_1,\dots,i_k)}(a))\otimes
x+(-1)^{s+k-1}\delta^{(i_1,\dots,i_k)}(a)\otimes d(x)=$$
$$=(-1)^{k-1}\delta^{(i_1,\dots,i_k)}(\p(a))\otimes
x+(-1)^{s+k-1}\delta^{(i_1,\dots,i_k)}(a)\otimes d(x)+$$
$$+\sum_{\sigma\in\Sigma_k}\sum_{I_{\sigma}} (-1)^{{\rm
sign}(\sigma)+k(m-1)}\delta^{(\widehat{\sigma(i_{m+1})},\dots,\widehat{\sigma(i_k)})}
(\delta^{(\widehat{\sigma(i_1)},\dots,\widehat{\sigma(i_m)})}(a))\otimes
x.$$ On the other hand since
$d\p_{(i_1,\dots,i_k)}+\p_{(i_1,\dots,i_k)}d=d(\p_{(i_1,\dots,i_k)})$,
by using the formula $(2.5)$ we obtain
$$d((-1)^{s(k-1)+(n+q)k}a\otimes\p_{(i_1,\dots,i_k)}(x))=$$
$$=(-1)^{s(k-1)+(n+q)k}\p(a)\otimes\p_{(i_1,\dots,i_k)}(x)+
(-1)^{sk+(n+q)k}a\otimes d(\p_{(i_1,\dots,i_k)}(x))=$$
$$=(-1)^{s(k-1)+(n+q)k}\p(a)\otimes\p_{(i_1,\dots,i_k)}(x)+
(-1)^{sk+(n+q)k+1}a\otimes\p_{(i_1,\dots,i_k)}(d(x))+$$
$$+(-1)^{sk+(n+q)k}a\otimes
\sum_{\sigma\in\Sigma_k}\sum_{I_{\sigma}} (-1)^{{\rm
sign}(\sigma)+1}
\p_{(\widehat{\sigma(i_1)},\dots,\widehat{\sigma(i_m)})}
(\p_{(\widehat{\sigma(i_{m+1})},\dots,\widehat{\sigma(i_k)})}(x)).$$
Now, by using the relations $(3.8)$ we obtain
$$(-1)^{k-1}\delta^{(i_1,\dots,i_k)}(\p(a))\otimes x\sim
(-1)^{k-1}(-1)^{(s-1)(k-1)+(n+q)k}\p(a)\otimes\p_{(i_1,\dots,i_k)}(x)=$$
$$=(-1)^{s(k-1)+(n+q)k}\p(a)\otimes\p_{(i_1,\dots,i_k)}(x),$$
$$(-1)^{s+k-1}\delta^{(i_1,\dots,i_k)}(a)\otimes d(x)\sim
(-1)^{s+k-1}(-1)^{s(k-1)+(n+q-1)k}a\otimes\p_{(i_1,\dots,i_k)}(d(x))=$$
$$=(-1)^{sk+(n+q)k+1}a\otimes\p_{(i_1,\dots,i_k)}(d(x)),$$
$$\sum_{\sigma\in\Sigma_k}\sum_{I_{\sigma}} (-1)^{{\rm
sign}(\sigma)+k(m-1)}\delta^{(\widehat{\sigma(i_{m+1})},\dots,\widehat{\sigma(i_k)})}
(\delta^{(\widehat{\sigma(i_1)},\dots,\widehat{\sigma(i_m)})}(a))\otimes
x\sim $$ $$\sim \sum_{\sigma\in\Sigma_k}\sum_{I_{\sigma}}
(-1)^{{\rm
sign}(\sigma)+k(m-1)}(-1)^\nu\delta^{(\widehat{\sigma(i_1)},\dots,\widehat{\sigma(i_m)})}(a)\otimes
\p_{(\widehat{\sigma(i_{m+1})},\dots,\widehat{\sigma(i_k)})}(x)\sim$$
$$\sim \sum_{\sigma\in\Sigma_k}\sum_{I_{\sigma}} (-1)^{{\rm
sign}(\sigma)+k(m-1)+\nu}(-1)^\varepsilon a\otimes
\p_{(\widehat{\sigma(i_1)},\dots,\widehat{\sigma(i_m)})}
(\p_{(\widehat{\sigma(i_{m+1})},\dots,\widehat{\sigma(i_k)})}(x)),$$
where $\nu=(s+m-1)(k-m-1)+(n+q)(k-m)$,
$\varepsilon=s(m-1)+(n+q-1)m$, and hence we have the relation
$$\sum_{\sigma\in\Sigma_k}\sum_{I_{\sigma}} (-1)^{{\rm
sign}(\sigma)+k(m-1)}\delta^{(\widehat{\sigma(i_{m+1})},\dots,\widehat{\sigma(i_k)})}
(\delta^{(\widehat{\sigma(i_1)},\dots,\widehat{\sigma(i_m)})}(a))\otimes
x\sim $$ $$\sim (-1)^{sk+(n+q)k}a\otimes
\sum_{\sigma\in\Sigma_k}\sum_{I_{\sigma}} (-1)^{{\rm
sign}(\sigma)+1}
\p_{(\widehat{\sigma(i_1)},\dots,\widehat{\sigma(i_m)})}
(\p_{(\widehat{\sigma(i_{m+1})},\dots,\widehat{\sigma(i_k)})}(x)),$$
since $\nu+\varepsilon\equiv sk+(n+q)k+k(m-1)+1\,{\rm mod}(2)$
holds. The verification of the required condition is
completed.~~~$\blacksquare$

{\bf Definition 3.1}. By a chain realization of the differential
module with $\infty$-sim\-plicial faces $(X,d,\widetilde{\p})$ we
mean the differential graded module $(|X|,d)$ that is defined as
the differential quotient module $$(|X|,d)=\bigoplus_{n\geqslant
0}(\overline{\F[n]}_\bu\otimes X_{n,\bu},d)/\!\sim$$ by the
equivalence relation $\sim\,$, which is generated by the relations
$(3.8)$.

Consider the differential module $(|X|,d)$ in more detail. Since
$\F[n]$ is a free $\F$-module with the generator $1_n$, from the
easily verified equality
$\delta^{(i_1,\dots,i_k)}(1_{n-k})=(-1)^{(n-k)k}\p_{(i_1,\dots,i_k)}(1_n)$
and the formulas $(3.4)$, $(3.8)$ follows that, for any elements
$a=\p_{(i_1,\dots,i_s)}\dots\p_{(j_1,\dots,
j_t)}(1_n)\in\overline{\F[n]}$ and $x\in X_{n,\bu}$, the relation
$a\otimes x\sim (-1)^\varepsilon 1_m\otimes
\p_{(i_1,\dots,i_s)}\dots\p_{(j_1,\dots, j_t)}(x)$ holds, where
$(-1)^\varepsilon$ is computed by using $(3.4)$ and $(3.8)$. In
particular, by using the equality $d(1_n)=0$ we obtain that, for
each element $x\in X_{n,q}$, the following relation holds:
$$\p(1_n)\otimes x=\sum_{0\leqslant i_1<\dots<i_k\leqslant
n}(-1)^{i_1+\dots+i_k}\p_{(i_1,\dots,i_k)}(1_n)\otimes x\sim$$
$$\sim \sum_{0\leqslant i_1<\dots<i_k\leqslant
n}(-1)^{i_1+\dots+i_k+n(k-1)+(q-1)k}
1_{n-k}\otimes\p_{(i_1,\dots,i_k)}(x).\eqno(3.9)$$ Let us consider
for $(X,d,\widetilde{\p})$ its corresponding the $\D$-module
$(X,d^k)$, and let us define new the $\D$-module
$(X,\overline{d}\,^k)$ by setting
$\overline{d}\,^k(x)=(-1)^{n(k-1)+(q-1)k}d^k(x)$, $k\geqslant 0$,
$x\in X_{n,q}$. A direct check shows that $(X,\overline{d}\,^k)$
really is a $\D$-module, and hence for $(X,\overline{d}\,^k)$ is
defined its corresponding differential module
$(\overline{X},\overline{\p})$, whose the differential
$\overline{\p}:\overline{X}_\bu\to\overline{X}_{\bu-1}$, for any
elements $x\in X_{n,q}\subset\overline{X}_{n+q}$, is given by the
formula $$\overline{\p}(x)=\sum_{k\geqslant
0}\,\overline{d}\,^k(x)=(-1)^nd(x)+\sum_{0\leqslant
i_1<\dots<i_k\leqslant n}(-1)^{i_1+\dots+i_k+n(k-1)+(q-1)k}
\p_{(i_1,\dots,i_k)}(x)$$ and extended by linearity to all
elements of the module $\overline{X}$.  In what follows the
elements $x\in X_{n,q}\subset\overline{X}_{n+q}$ we will denote by
$[x]$. We note that each equivalence class $[a\otimes x]\in|X|$ is
a linear combination of equivalence classes of the form
$[1_n\otimes y]$, where the elements $y\in X_{n,\bu}$ are uniquely
determined by the relations $(3.4)$ and $(3.8)$. Consider the map
of graded modules $F:|X|_\bu\to \overline{X}_\bu$ that, for any
elements $[1_n\otimes x]$, where $x\in X_{n,q}$, is given by the
equality $F([1_n\otimes x])=[x]$ and extended by linearity to all
elements of the module $|X|$. Clear that the map $F$ is an
isomorphism of graded modules. Furthermore, from the relations
$(3.9)$ follows that $F$ is an isomorphism of differential modules
$F:(|X|,d)\to (\overline{X},\overline{\p})$.

Now we show that each $\A$-algebra defines some the differential
module with $\infty$-simplicial faces $(T(A),d,\widetilde{\p})$.
Recall \cite{S} (see also \cite{Kad}, \cite{Smirnov}) that an
$\A$-algebra $(A,d,\pi_n)$ is defined as a differential module
$(A,d)$, $A=\{A_n\}$, $n\in\mathbb{Z}$, $n\geqslant 0$,
$d:A_\bu\to A_{\bu-1}$, equipped with a family of maps
$\{\pi_n:(A^{\otimes (n+2)})_\bu\to
A_{\bu+n}\,\,|\,\,n\in\mathbb{Z},\,\,n\geqslant 0\}$ which, for
all integers $n\geqslant -1$, satisfy the following relations:
$$d(\pi_{n+1})=\sum\limits_{m=0}^n\sum_{t=1}^{m+2}(-1)^{t(n-m+1)+n+1}
\pi_m(\underbrace{1\otimes\dots\otimes
1}_{t-1}\otimes\,\pi_{n-m}\otimes\underbrace{1\otimes\dots \otimes
1}_{m-t+2}),\eqno(3.10)$$ where
$d(\pi_{n+1})=d\pi_{n+1}+(-1)^n\pi_{n+1}d$. For any $\A$-algebra
$(A,d,\pi_n)$, we consider the differential bigraded module
$(T(A),d)$, where $T(X)=\{T(A)_{n,m}\}$, $n,m\in\mathbb{Z}$,
$n\geqslant 0$, $m\geqslant 0$, which is defined by the equalities
$T(A)_{n,m}=(A^{\otimes n})_m$, $n>0$, $m>0$, $T(A)_{0,0}=K$,
$T(A)_{n,0}=0$, $n>0$, $T(A)_{0,m}=0$, $m>0$, and its the
differential $d:T(A)_{n,\bu}\to T(A)_{n,\bu-1}$ is the usual
differential in a tensor product. We define the family of the maps
$\widetilde{\p}=\{\p^n_{(i_1,\dots,i_k)}:T(A)_{n,q}\to
T(A)_{n-k,q+k-1}\}$, where $n\geqslant 0$, $q\geqslant 0$ and
$0\leqslant i_1<\dots<i_k\leqslant n$, by setting
$$\p^n_{(i_1,\dots,i_k)}\!= \left\{\begin{array}{ll}
(-1)^{k(q-1)}1^{\otimes (j-1)}\otimes\,\pi_{k-1}\otimes 1^{\otimes
(n-k-j)},\quad {\rm{if}}\quad 1\leqslant j\leqslant n-k&\\
{\rm{and}}\quad (i_1,\dots,i_k)=(j,j+1,\dots,j+k-1);\\ 0,\quad
{\rm{otherwise}}.&\\
\end{array}\right.\eqno(3.11)$$

{\bf Theorem 3.2}. Given any $\A$-algebra $(A,d,\pi_n)$, the
triple $(T(A),d,\widetilde{\p})$ is a differential module with
$\infty$-simplicial faces.

{\bf Proof}. For a family of maps
$\widetilde{\p}=\{\p^n_{(i_1,\dots,i_k)}:T(A)_{n,q}\to
T(A)_{n-k,q+k-1}\}$, which are defined by the formula $(3.11)$, we
need to check the validity of the relations $(2.5)$. First, for
the maps
$\p^{n+3}_{(1,2,\dots,n+2)}=(-1)^{(n+2)(q-1)}\pi_{n+1}:(A^{\otimes(n+3})_q\to
A_{q+n+1}$, $n\geqslant -1$, we will verify that the relations
$(2.5)$ are true. Clearly, for $n=-1$, i.e. for
$\p^2_{(1)}=(-1)^{q-1}\pi_0$, the condition $(2.5)$ holds, since
from the equality $d(\pi_0)=d\pi_0-\pi_0d=0$ we obtain
$d(\p^2_{(1)})=d\p^2_{(1)}+\p^2_{(1)}d=(-1)^{q-1}d(\pi_0)=0$.
Suppose now that $n\geqslant 0$. We note that, for any permutation
$\sigma\in\Sigma_{n+2}$ of integers $1,\dots,n+2$, the equality
$\widehat{\sigma(1)}=1$ is always true. Indeed, for the collection
$(\sigma(1),\dots,\sigma(n+2))$, all its numbers smaller than
$\sigma(1)$ are to the right of $\sigma(1)$, and the quantity of
such numbers is $\alpha(\sigma(1))=\sigma(1)-1$. From the
condition $\widehat{\sigma(1)}=1$ and the formula $(3.11)$ follows
that in considered case the relations $(2.5)$ can be written as
$$d(\p^{n+3}_{(1,2,\dots,n+2)})=\sum_{m=0}^n\sum_{t=1}^{m+2}(-1)^{{\rm
sign}(\sigma_{t,m})+1}\p^{m+2}_{(1,2,\dots,m+1)}\p^{n+3}_{(t,t+1,\dots,t+n-m)},\eqno(3.12)$$
where
$d(\p^{n+3}_{(1,2,\dots,n+2)})=d\p^{n+3}_{(1,2,\dots,n+2)}+\p^{n+3}_{(1,2,\dots,n+2)}d$,
and $\sigma_{t,m}\in\Sigma_{n+2}$ is a permutation of numbers
$1,\dots,n+2$, which satisfies the conditions
$$\widehat{\sigma_{t,m}(1)}=1,\,\widehat{\sigma_{t,m}(2)}=
2,\dots,\,\widehat{\sigma_{t,m}(m+1)}=m+1,$$
$$\widehat{\sigma_{t,m}(m+2)}=t,\,\widehat{\sigma_{t,m}(m+3)}=
t+1,\dots,\widehat{\sigma_{t,m}(n+2)}=t+n-m.$$ It is easy to
verify that the permutation $\sigma_{t,m}$ acts on the collection
$(1,2,\dots,n+2)$ by a partitioning of this collection on three
blocks
$$(1,2,\dots,t-1\,|\,t,t+1,\dots,t+n-m\,|\,t+n-m+1,\dots,n+2)$$
and by a permutation of places of the second and third blocks.
From this follows that the number of inversions $I(\sigma_{t,m})$
of the permutation $\sigma_{t,m}$ is equal to the product of the
lengths of the second and third blocks of the above partition,
i.e. $I(\sigma_{t,m})=(n-m+1)(m-t+2)$. If now multiply the left
and right sides of the equality $(3.10)$ by $(-1)^{(n+2)(q-1)}$
and notice that $I(\sigma_{t,m})\equiv {\rm
sign}(\sigma_{t,m})\,{\rm mod}(2)$, then by using $(3.11)$ we
obtain the relations $(3.12)$. In an analogous way we can verify
the validity of the relations $(2.5)$ for an arbitrary maps
$\p^n_{(j,j+1,\dots,j+k-1)}$, where $k\geqslant 1$ and $1\leqslant
j\leqslant n-k$. Consider now the case of the maps
$\p^n_{(i_1,\dots,i_k)}$, when the ordered collection
$(i_1,\dots,i_k)$ is not a collection of the form
$(j,j+1,\dots,j+k-1)$, where $k\geqslant 1$ and $1\leqslant
j\leqslant n-k$. It follows from the Proposition 2.1 that in this
case the relations $(2.5)$ require a verification only for those
collections $(i,\dots,i_k)$, which have the form
$$(j,j+1,\dots,j+l-1,t,t+1,\dots,t+m-1),$$ where $l\geqslant 1$,
$m\geqslant 1$, $l+m=k$, $j+l<t$, since in otherwise the left and
right sides of the relations $(2.5)$ are equal zero. For
collections of the above form, the relations $(2.5)$ can be
written as the relations
$$d(\p^n_{(j,j+1,\dots,j+l-1,t,t+1,\dots,t+m-1)})=0=
-\p^{n-m}_{(j,j+1,\dots,j+l-1)}\p^n_{(t,t+1,\dots,t+m-1)}+$$
$$+(-1)^{lm+1}\p^{n-l}_{(t-l,t+1-l,\dots,t+m-1-l)}\p^n_{(j,j+1,\dots,j+l-1)},$$
that easy follows from the obvious equalities
$$(\underbrace{1\otimes\dots\otimes
1}_{j-1}\otimes\,\pi_{l-1}\otimes\underbrace{1\otimes\dots\otimes
1}_{n-m-j-l})(\underbrace{1\otimes\dots\otimes
1}_{t-1}\otimes\,\pi_{m-1}\otimes\underbrace{1\otimes\dots\otimes
1}_{n-m-t})=$$
$$=(-1)^{(m-1)(l-1)}(\underbrace{1\otimes\dots\otimes
1}_{t-l-1}\otimes\,\pi_{m-1}\otimes\underbrace{1\otimes\dots\otimes
1}_{n-m-t})(\underbrace{1\otimes\dots\otimes
1}_{j-1}\otimes\,\pi_{l-1}\otimes\underbrace{1\otimes\dots\otimes
1}_{n-l-j}),$$ where $l\geqslant 1$, $m\geqslant 1$,
$j+l<t$.~~~$\blacksquare$

For any $\A$-algebra $(A,d,\pi_n)$, let us consider the chain
realization $(|T(A)|,d)$ of the differential module with
$\infty$-simplicial faces $(T(A),d,\widetilde{\p})$, which
corresponds to the $\A$-algebra $(A,d,\pi_n)$. As was said above,
by using the isomorphism $F$ we can identify the differential
module $(|T(A)|,d)$ and the differential module
$(\overline{T(A)},\overline{\p})$; here
$\overline{T(A)}_m=\bigoplus_{n+q=m}(A^{\otimes n})_q$ and the
differential $\overline{\p}:\overline{T(A)}_\bu\to
\overline{T(A)}_{\bu-1}$ is given by the following formula:
$$\overline{\p}([a_1\otimes\dots\otimes
a_n])=(-1)^n\sum_{i=1}^n(-1)^{\varepsilon}[a_1\otimes\dots\otimes
d(a_i)\otimes\dots\otimes a_n]\,+\eqno(3.13)$$
$$+\sum_{k=1}^{n-1}\sum_{i=1}^{n-k}(-1)^{\frac{k(k-1)}{2}+ik+n(k-1)+\varepsilon(k-1)}
[a_1\otimes\dots\otimes \pi_{k-1}(a_i\otimes\dots\otimes
a_{i+k})\otimes\dots\otimes a_n],$$ where
$\varepsilon=\deg(a_1)+\dots+\deg(a_{i-1})$ and elements
$[a_1\otimes\dots\otimes a_n]\in\overline{T(A)}_{n+q}$ are
generators of the module $\overline{T(A)}$.

Let us proceed to defining the structure of a differential
coalgebra $(|T(A)|,d,\nabla)$ on the differential module
$(|T(A)|,d)$, where $(A,d,\pi_n)$ is an arbitrary $\A$-algebra.
For this we introduce a conception of the tensor product of
differential modules with $\infty$-simplicial faces.

{\bf Definition 3.2}. By a tensor product $(X\otimes
Y,d,\widetilde{\p})$ of the differential modules with
$\infty$-simplicial faces $(X,d,\widetilde{\p})$ and
$(Y,d,\widetilde{\p})$ we mean a tensor product $(X\otimes Y,d)$
of the differential bigraded modules $(X,d)$ and $(Y,d)$ equipped
with a family of the $\infty$-simplicial faces
$\widetilde{\p}=\{\p_{(i_1,\dots,i_k)}:(X\otimes
Y)_{n,\bu}\to(X\otimes Y)_{n-k,\bu+k-1}\}$, which are given on an
arbitrary element $x\otimes y\in X_{q,s}\otimes Y_{l,t}$ by the
following rule: $$\p_{(i_1,\dots,i_k)}(x\otimes y)=\eqno(3.14)$$
$$=\left\{\begin{array}{ll} \p_{(i_1,\dots,i_k)}(x)\otimes
y,&0\leqslant i_1<\dots<i_k<q,\\ 0,&0\leqslant i_1<\dots<i_k=q,\\
(-1)^{(k-1)q+s}x\otimes\p_{(i_1-q,\dots,i_k-q)}(y),&q<i_1<\dots<i_k\leqslant
q+l,\\ 0,&0\leqslant i_1\leqslant q<i_k\leqslant q+l.\\
\end{array}\right.$$

{\bf Theorem 3.3}. A tensor product of differential modules with
$\infty$-simplicial faces is a differential module with
$\infty$-simplicial faces.

{\bf Proof}. For $\infty$-simplicial faces in a tensor product
$(X\otimes Y ,d,\widetilde{\p})$ of the differential modules with
$\infty$-simplicial faces $(X,d,\widetilde{\p})$ and
$(Y,d,\widetilde{\p})$, we need to check the validity of the
relations $(2.5)$. We will check each case of the formula
$(3.14)$.

1). Suppose that $0\leqslant i_1<\dots<i_k<q$. In this case, for
any element $x\otimes y\in X_{q,s}\otimes Y_{l,t}$, we obtain
$$(d(\p_{(i_1,\dots,i_k)}))(x\otimes y)=$$
$$=d(\p_{(i_1,\dots,i_k)}(x\otimes
y))+(-1)^{-k+(k-1)+1}\p_{(i_1,\dots,i_k)}(d(x\otimes y))=$$
$$=d(\p_{(i_1,\dots,i_k)}(x)\otimes
y)+\p_{(i_1,\dots,i_k)}(d(x)\otimes y+(-1)^{q+s}x\otimes d(y))=$$
$$=d(\p_{(i_1,\dots,i_k)}(x))\otimes
y+(-1)^{q+s-1}\p_{(i_1,\dots,i_k)}(x)\otimes d(y)+$$
$$+\,\p_{(i_1,\dots,i_k)}(d(x))\otimes
y+(-1)^{q+s}\p_{(i_1,\dots,i_k)}(x)\otimes
d(y)=(d(\p_{(i_1,\dots,i_k)}))(x)\otimes y=$$
$$=\sum_{\sigma\in\Sigma_k}\sum_{I_{\sigma}} (-1)^{{\rm
sign}(\sigma)+1}
\p_{(\widehat{\sigma(i_1)},\dots,\widehat{\sigma(i_m)})}
\p_{(\widehat{\sigma(i_{m+1})},\dots,\widehat{\sigma(i_k)})}(x)\otimes
y=$$ $$=\sum_{\sigma\in\Sigma_k}\sum_{I_{\sigma}} (-1)^{{\rm
sign}(\sigma)+1}
(\p_{(\widehat{\sigma(i_1)},\dots,\widehat{\sigma(i_m)})}\otimes
1)(\p_{(\widehat{\sigma(i_{m+1})},\dots,\widehat{\sigma(i_k)})}\otimes
1)(x\otimes y)=$$ $$=\sum_{\sigma\in\Sigma_k}\sum_{I_{\sigma}}
(-1)^{{\rm sign}(\sigma)+1}
\p_{(\widehat{\sigma(i_1)},\dots,\widehat{\sigma(i_m)})}
\p_{(\widehat{\sigma(i_{m+1})},\dots,\widehat{\sigma(i_k)})}(x\otimes
y).$$

2). Suppose that $0\leqslant i_1<\dots<i_k=q$. In this case, for
an arbitrary element $x\otimes y\in X_{q,s}\otimes Y_{l,t}$, by
the definition of $\infty$-simplicial faces in a tensor product we
have $$(d(\p_{(i_1,\dots,i_k)}))(x\otimes y)=0.$$ Now we need
check that if $0\leqslant i_1<\dots<i_k=q$, then for any element
$x\otimes y\in X_{q,s}\otimes Y_{l,t}$ the following equality
holds: $$\sum_{\sigma\in\Sigma_k}\sum_{I_{\sigma}} (-1)^{{\rm
sign}(\sigma)+1}
\p_{(\widehat{\sigma(i_1)},\dots,\widehat{\sigma(i_m)})}
\p_{(\widehat{\sigma(i_{m+1})},\dots,\widehat{\sigma(i_k)})}(x\otimes
y)=0.\eqno(3.15)$$ For each summand of the left side of $(3.15)$,
we will show that the following equality holds:
$$\p_{(\widehat{\sigma(i_1)},\dots,\widehat{\sigma(i_m)})}
\p_{(\widehat{\sigma(i_{m+1})},\dots,\widehat{\sigma(i_k)})}(x\otimes
y)=0.$$ Indeed, for an arbitrary summand
$\p_{(\widehat{\sigma(i_1)},\dots,\widehat{\sigma(i_m)})}
\p_{(\widehat{\sigma(i_{m+1})},\dots,\widehat{\sigma(i_k)})}(x\otimes
y)$, we suppose that $\sigma\in\Sigma_k$ is a permutation such
that $i_k=q\in\{\sigma(i_{m+1}),\dots,\sigma(i_k)\}$. In this case
we will show that $\widehat{\sigma(i_k)}=q$ is true. Assume by
contradiction that we have $\widehat{\sigma(i_k)}<q$. From this,
by using $\widehat{\sigma(i_k)}=\sigma(i_k)$, we obtain
$\sigma(i_k)<q$. Since for the $\infty$-simplicial face
$\p_{(\widehat{\sigma(i_{m+1})},\dots,\widehat{\sigma(i_k)})}$ the
condition $0\leqslant\widehat{\sigma(i_{m+1})}<\dots
<\widehat{\sigma(i_k)}<q$ holds, we get
$\widehat{\sigma(i_{k-p})}<q-p$, $0\leqslant p\leqslant k-(m+1)$.
Since
$\widehat{\sigma(i_{k-p})}=\sigma(i_{k-p})-\alpha(\sigma(i_{k-p}))$
and $0\leqslant\alpha(\sigma(i_{k-p}))\leqslant p$, we obtain
$\sigma(i_{k-p})<q$, $0\leqslant p\leqslant k-(m+1)$, that
contradicts the condition
$q\in\{\sigma(i_{m+1}),\dots,\sigma(i_k)\}$. Thus, the above
assumption $\widehat{\sigma(i_k)}<q$ is incorrect and hence we
have $\widehat{\sigma(i_k)}=q$. From this follows that in
considered case, for each summand
$\p_{(\widehat{\sigma(i_1)},\dots,\widehat{\sigma(i_m)})}
\p_{(\widehat{\sigma(i_{m+1})},\dots,\widehat{\sigma(i_k)})}(x\otimes
y)$, where $x\otimes y\in X_{q,s}\otimes Y_{l,t}$, we obtain, by
using $\widehat{\sigma(i_k)}=q$, that the equality
$\p_{(\widehat{\sigma(i_1)},\dots,\widehat{\sigma(i_m)})}
\p_{(\widehat{\sigma(i_{m+1})},\dots,\widehat{\sigma(i_k)})}(x\otimes
y)=0$ is true. Now, for an arbitrary summand
$\p_{(\widehat{\sigma(i_1)},\dots,\widehat{\sigma(i_m)})}
\p_{(\widehat{\sigma(i_{m+1})},\dots,\widehat{\sigma(i_k)})}(x\otimes
y)$, we suppose that $\sigma\in\Sigma_k$ is a permutation such
that $i_k=q\in\{\sigma(i_1),\dots,\sigma(i_m)\}$. In this case we
will show that $\widehat{\sigma(i_m)}=q-(k-m)$ is true. Assume by
contradiction that we have $\widehat{\sigma(i_m)}<q-(k-m)$. From
this, by using
$\widehat{\sigma(i_m)}=\sigma(i_m)-\alpha(\sigma(i_m))$ and
$0\leqslant\alpha(\sigma(i_m))\leqslant k-m$, we obtain
$\sigma(i_m)<q$. Since for the $\infty$-simplicial face
$\p_{(\widehat{\sigma(i_1)},\dots,\widehat{\sigma(i_m)})}$ the
condition $0\leqslant\widehat{\sigma(i_1)}<\dots
<\widehat{\sigma(i_m)}<q$ holds, we get
$\widehat{\sigma(i_p)}<q-(k-p)$, $1\leqslant p\leqslant m$. Since
$\widehat{\sigma(i_p)}=\sigma(i_p)-\alpha(\sigma(i_p))$ and
$0\leqslant\alpha(\sigma(i_p))\leqslant k-p$, we obtain
$\sigma(i_p)<q$, $1\leqslant p\leqslant m$, that contradicts the
condition $q\in\{\sigma(i_1),\dots,\sigma(i_m)\}$. Thus, the above
assumption $\widehat{\sigma(i_m)}<q-(k-m)$ is incorrect and hence
we have $\widehat{\sigma(i_m)}=q-(k-m)$. From this follows that in
considered case, for each summand
$\p_{(\widehat{\sigma(i_1)},\dots,\widehat{\sigma(i_m)})}
\p_{(\widehat{\sigma(i_{m+1})},\dots,\widehat{\sigma(i_k)})}(x\otimes
y)$, where $x\otimes y\in X_{q,s}\otimes Y_{l,t}$, we obtain, by
using $\widehat{\sigma(i_m)}=q-(k-m)$ and
$\p_{(\widehat{\sigma(i_{m+1})},\dots,\widehat{\sigma(i_k)})}(x\otimes
y)=\p_{(\widehat{\sigma(i_{m+1})},\dots,\widehat{\sigma(i_k)})}(x)\otimes
y\in X_{q-(k-m),s+(k-m)-1}\otimes Y_{l,t}$, that the equality
$\p_{(\widehat{\sigma(i_1)},\dots,\widehat{\sigma(i_m)})}
\p_{(\widehat{\sigma(i_{m+1})},\dots,\widehat{\sigma(i_k)})}(x\otimes
y)=0$ is true.

3). Suppose that $q<i_1<\dots<i_k\leqslant q+l$. In this case, for
an arbitrary element $x\otimes y\in X_{q,s}\otimes Y_{l,t}$, we
have  $$(d(\p_{(i_1,\dots,i_k)}))(x\otimes y)=$$
$$=d(\p_{(i_1,\dots,i_k)}(x\otimes
y))+(-1)^{-k+(k-1)+1}\p_{(i_1,\dots,i_k)}(d(x\otimes y))=$$
$$d((-1)^{q(k-1)+s}x\otimes
\p_{(i_1-q,\dots,i_k-q)}(y))+\p_{(i_1,\dots,i_k)}(d(x)\otimes
y+(-1)^{q+s}x\otimes d(y))=$$ $$=(-1)^{q(k-1)+s} d(x)\otimes
\p_{(i_1-q,\dots,i_k-q)}(y)+(-1)^{qk}x\otimes
d(\p_{(i_1-q,\dots,i_k-q)}(y))+$$
$$+\,(-1)^{q(k-1)+s-1}d(x)\otimes
\p_{(i_1-q,\dots,i_k-q)}(y)+(-1)^{qk}x\otimes
\p_{(i_1-q,\dots,i_k-q)}(d(y))=$$ $$=(-1)^{qk}x\otimes
(d(\p_{(i_1-q,\dots,i_k-q)}))(y)=$$
$$=(-1)^{qk}\sum_{\sigma\in\Sigma_k}\sum_{I_{\sigma}} (-1)^{{\rm
sign}(\sigma)+1}x\otimes
\p_{(\widehat{\sigma(i_1-q)},\dots,\widehat{\sigma(i_m-q)})}
\p_{(\widehat{\sigma(i_{m+1}-q)},\dots,\widehat{\sigma(i_k-q)})}(y)=$$
$$=(-1)^{qk}\sum_{\sigma\in\Sigma_k}\sum_{I_{\sigma}} (-1)^{{\rm
sign}(\sigma)+1}x\otimes
\p_{(\widehat{\sigma(i_1)}-q,\dots,\widehat{\sigma(i_m)}-q)}
\p_{(\widehat{\sigma(i_{m+1})}-q,\dots,\widehat{\sigma(i_k)}-q)}(y)=$$
$$=\sum_{\sigma\in\Sigma_k}\sum_{I_{\sigma}} (-1)^{{\rm
sign}(\sigma)+1}
\p_{(\widehat{\sigma(i_1)},\dots,\widehat{\sigma(i_m)})}
\p_{(\widehat{\sigma(i_{m+1})},\dots,\widehat{\sigma(i_k)})}(x\otimes
y).$$

4). Suppose that $0\leqslant i_1\leqslant q<i_k$. In this case,
for any element $x\otimes y\in X_{q,s}\otimes Y_{l,t}$, by the
definition of $\infty$-simplicial faces in a tensor product, we
have $$(d(\p_{(i_1,\dots,i_k)}))(x\otimes y)=0.$$ Now we need
check that if $0\leqslant i_1\leqslant q<i_k$, then for any
element $x\otimes y\in X_{q,s}\otimes Y_{l,t}$ the equality
$(3.15)$ is true. For any summand $(-1)^{{\rm
sign}(\sigma)+1}\p_{(\widehat{\sigma(i_1)},\dots,\widehat{\sigma(i_m)})}
\p_{(\widehat{\sigma(i_{m+1})},\dots,\widehat{\sigma(i_k)})}(x\otimes
y)$ in $(3.15)$, we suppose that $\sigma\in\Sigma_k$ is a
permutation such that
$i_k=q+c\in\{\sigma(i_{m+1}),\dots,\sigma(i_k)\}$, $c=i_k-q>0$. In
this case, similar to how it was done in 2), we can show that
$\widehat{\sigma(i_k)}=q+c$, i.e. that $\widehat{\sigma(i_k)}>q$.
If $\widehat{\sigma(i_{m+1})}\leqslant q$, then $(-1)^{{\rm
sign}(\sigma)+1}\p_{(\widehat{\sigma(i_1)},\dots,\widehat{\sigma(i_m)})}
\p_{(\widehat{\sigma(i_{m+1})},\dots,\widehat{\sigma(i_k)})}(x\otimes
y)=0$, since
$\p_{(\widehat{\sigma(i_{m+1})},\dots,\widehat{\sigma(i_k)})}(x\otimes
y)=0$. If $\widehat{\sigma(i_{m+1})}>q$, then consider separately
the cases a)
$\widehat{\sigma(i_{m})}\geqslant\widehat{\sigma(i_{m+1})}$ and b)
$\widehat{\sigma(i_{m})}<\widehat{\sigma(i_{m+1})}$. In the case
a) we have $\widehat{\sigma(i_{m})}>q$.  If we assume that
$i_1=\sigma(i_l)$ for some $1\leqslant l\leqslant k$, then obtain
$\widehat{\sigma(i_l)}=\sigma(i_l)$. From the inequalities
$\widehat{\sigma(i_l)}=i_1\leqslant q$ and
$q<\widehat{\sigma(i_{m+1})}<\dots<\widehat{\sigma(i_k)}$ follows
that $1\leqslant l\leqslant m$ and hence we have
$\widehat{\sigma(i_1)}<\widehat{\sigma(i_l)}\leqslant q$. Thus, in
the case a) we obtain $$(-1)^{{\rm
sign}(\sigma)+1}\p_{(\widehat{\sigma(i_1)},\dots,\widehat{\sigma(i_m)})}
\p_{(\widehat{\sigma(i_{m+1})},\dots,\widehat{\sigma(i_k)})}(x\otimes
y)=$$ $$=(-1)^{{\rm
sign}(\sigma)+1+(k-m-1)q+s}\p_{(\widehat{\sigma(i_1)},\dots,\widehat{\sigma(i_m)})}(x\otimes
\p_{(\widehat{\sigma(i_{m+1})}-q,\dots,\widehat{\sigma(i_k)}-q)}(y))=0.$$
In the case b) we have
$\widehat{\sigma(i_1)}<\dots<\widehat{\sigma(i_m)}<\widehat{\sigma(i_{m+1})}<\dots<\widehat{\sigma(i_k)}$.
From this inequalities, similar to how it was done in the proof of
Proposition 2.1, we obtain the inequalities
$\sigma(i_1)<\dots<\sigma(i_m)<\sigma(i_{m+1})<\dots<\sigma(i_k)$,
which imply that $\sigma$ is the identity permutation. If $i_m>q$,
then we get
$(-1)^1\p_{(i_1,\dots,i_m)}\p_{(i_{m+1},\dots,i_k)}(x\otimes
y)=0$, since $i_1\leqslant q$. If $i_m<q$, then in the sum
$(3.15)$ we obtain the summand
$$(-1)^1\p_{(i_1,\dots,i_m)}\p_{(i_{m+1},\dots,i_k)}(x\otimes
y)=(-1)^{(k-m-1)q+s+1}\p_{(i_1,\dots,i_m)}\otimes\p_{(i_{m+1}-q,\dots,i_k-q)}.$$
In $(3.15)$ this summand vanishes together with $(-1)^{{\rm
sign}(\varrho)+1}\p_{(i_{m+1}-m,\dots,i_k-m)}\p_{(i_1,\dots,i_m)}$,
where $\varrho$ is a permutation, which permute the places the
blocks $(i_1,\dots,i_m)$ and $(i_{m+1},\dots,i_k)$. Indeed, since
${\rm sign}(\varrho)=m(k-m)$, we obtain $$(-1)^{{\rm
sign}(\varrho)+1}\p_{(i_{m+1}-m,\dots,i_k-m)}\p_{(i_1,\dots,i_m)}(x\otimes
y)=$$
$$=(-1)^{m(k-m)+1+(k-m-1)(q-m)+s+m-1}\p_{(i_1,\dots,i_m)}\otimes\p_{(i_{m+1}-q,\dots,i_k-q)}.$$
Thus, for $(-1)^{{\rm
sign}(\sigma)+1}\p_{(\widehat{\sigma(i_1)},\dots,\widehat{\sigma(i_m)})}
\p_{(\widehat{\sigma(i_{m+1})},\dots,\widehat{\sigma(i_k)})}(x\otimes
y)$ in the sum $(3.15)$, if the permutation $\sigma\in\Sigma_k$
such that $i_k=q+c\in\{\sigma(i_{m+1}),\dots,\sigma(i_k)\}$,
$c=i_k-q\geqslant 0$, then this summand either is equal to zero or
vanishes together with another summand in this sum.  In a similar
way we verify that, for $(-1)^{{\rm
sign}(\sigma)+1}\p_{(\widehat{\sigma(i_1)},\dots,\widehat{\sigma(i_m)})}
\p_{(\widehat{\sigma(i_{m+1})},\dots,\widehat{\sigma(i_k)})}(x\otimes
y)$ in the sum $(3.15)$, if $\sigma\in\Sigma_k$ such that
$i_k=q+c\in\{\sigma(i_1),\dots,\sigma(i_m)\}$, $c=i_k-q\geqslant
0$, then this summand either is equal to zero or vanishes together
with another summand in this sum.

Thus, in all cases, for $\infty$-simplicial faces of the tensor
product $(X\otimes Y ,d,\widetilde{\p})$, the condition $(2.5)$
holds, and hence $(X\otimes Y,d,\widetilde{\p})$ is a differential
module with $\infty$-simplicial faces.~~~$\blacksquare$

Now, by using the concept of a tensor product of differential
modules with $\infty$-simplicial faces, we prove the following
assertion.

{\bf Theorem 3.4}. For any $\A$-algebra $(A,d,\pi_n)$, the chain
realization $(|T(A)|,d)$ of the differential $\F$-module
$(T(A),d,\widetilde{\p})$ equipped with a comultiplication
$\nabla:|T(A)|\to |T(A)|\otimes |T(A)|$, which on generators is
given by the formula $$\nabla([1_n\otimes (a_1\otimes\dots\otimes
a_n)])=$$ $$=\sum_{i=0}^{n}(-1)^{(n-i)\varepsilon_i}[1_i\otimes
(a_1\otimes\dots\otimes a_i)]\otimes [1_{n-i}\otimes
(a_{i+1}\otimes\dots\otimes a_n)],$$ where
$\varepsilon_i=\deg(a_1)+\dots+\deg(a_i)$, is a differential
coalgebra.

{\bf Proof}. It is easy to see that the comultiplication $\nabla$
is an associative. Let us show that $\nabla$ is a map of
differential modules. For the differential module with
$\infty$-simplicial faces $(T(A),d,\widetilde{\p})$, we consider
the differential module with $\infty$-simplicial faces
$(T(A),d\,',\widetilde{\p\,'})$, which is defined by the
equalities $d\,'=(-1)^n d:T(A)_{n,\bu}\to T(A)_{n,\bu-1}$ and
$$\widetilde{\p}\,'=\{\p\,'_{(i_1,\dots,i_k)}=(-1)^{n(k-1)+k(q-1)}\p_{(i_1,\dots,i_k)}:T(A)_{n,q}\to
T(A)_{n-k,q+k-1}\}.$$ Similarly to the proof of the Theorem 3.2
can be proved that the triple $(T(A),d\,',\widetilde{\p\,'})$ is a
differential module with $\infty$-simplicial faces. Direct
calculations show that the family of maps
$$\widetilde{\Delta}=\{\Delta_{(i_1,\dots,i_k)}:T(A)_{n,\bu}\to
(T(A)\otimes T(A))_{n-k,\bu+k}\},\quad
\Delta_{(i_1,\dots,i_k)}=0,~k\geqslant 1,$$
$$\Delta_{(\,\,)}(a_1\otimes\dots\otimes
a_n)=\sum_{i=0}^n(-1)^{(n-i)\varepsilon_i}(a_1\otimes\dots
a_i)\otimes (a_{i+1}\otimes\dots\otimes a_n),$$ where
$\varepsilon_i=\deg(a_1)+\dots+\deg(a_i)$, is the morphism of
differential modules with $\infty$-simplicial faces
$$\widetilde{\Delta}:(T(A),d\,',\widetilde{\p\,'})\to
(T(A),d\,',\widetilde{\p\,'})\otimes
(T(A),d\,',\widetilde{\p\,'})=(T(A)\otimes
T(A),d\,',\widetilde{\p\,'}).$$ From this it follows that the map
of bigraded modules $\Delta_{(\,\,)}$ induces the map of
differential modules $\Delta:(\overline{T(A)},\p\,')\to
(\overline{T(A)\otimes T(A)},\p\,')$. Now, by using the formula
$(3.11)$, we note that the $\infty$-simplicial faces
$\p\,'_{(i_1,\dots,i_k)}:T(A)_{n,\bu}\to T(A)_{n-k,\bu+k-1}$,
where $i_1=0$ or $i_k=n$, satisfy the condition
$\p\,'_{(i_1,\dots,i_k)}=0$. From this it follows that the
equality of differential modules $(\overline{T(A)\otimes
T(A)},\p\,')=(\overline{T(A)},\p\,')\otimes
(\overline{T(A)},\p\,')$ holds. By using this equality, and also
by using the equality
$(\overline{T(A)},\p\,')=(\overline{T(A)},\overline{\p})$, we
obtain the map of differential modules
$\Delta:(\overline{T(A)},\overline{\p})\to
(\overline{T(A)},\overline{\p})\otimes
\overline{T(A)},\overline{\p})$, which on generators is given by
the following formula: $$\Delta([a_1\otimes\dots\otimes
a_n])=\sum_{i=0}^n(-1)^{(n-i)\varepsilon_i}[a_1\otimes\dots\otimes
a_i]\otimes [a_{i+1}\otimes\dots\otimes a_n].\eqno(3.16)$$ Clear
that $\Delta$ satisfies the condition $(\Delta\otimes
1)\Delta=(1\otimes\Delta)\Delta$. Now, if we consider the
isomorphism of differential modules $F:(|T(A)|,d)\to
(\overline{T(A)},\overline{\p})$, $F([1_n\otimes x])=[x]$, then it
is easy to see that for the map $\nabla$, which is given in the
statement of this theorem, the equality $\nabla=(F^{-1}\otimes
F^{-1})\Delta F$ holds, and hence $\nabla$ is a map of
differential modules. Moreover, from this equality follows that
$F$ is an isomorphism of differential coalgebras
$F:(|T(A)|,d,\nabla)\to (\overline{T(A)},\overline{\p},\Delta)$.
~~~$\blacksquare$

The differential coalgebra
$(\overline{T(A)},\overline{\p},\Delta)$ is well known \cite{Kad}
as the $B$-construction $(B(A),d,\Delta)$ of the $\A$-algebra
$(A,d,\pi_n)$. Thus, up to isomorphism of differential coalgebras,
the $B$-construction of the $\A$-algebra $(A,d,\pi_n)$ is a chain
realization of the tensor differential module with
$\infty$-simplicial faces $(T(A),d,\widetilde{\p})$, which is
defined by this $\A$-algebra $(A,d,\pi_n)$.

We conclude this paragraph by pointing out that if an $\A$-algebra
$(A,d,\pi_n)$ is the differential associative algebra $(A,d,\pi)$,
where $\pi_0=\pi$ and $\pi_n=0$, $n>0$, then
$(\overline{T(A)},\overline{\p},\Delta)$ is usual $B$-construction
$(B(A),d,\Delta)$ of the algebra $(A,d,\pi)$, and for this
$B$-construction, the signs in the formulas $(3.13)$ and $(3.16)$
are coincide with the sings, which was given in \cite{Smirnov}.

\end{document}